\documentclass[12pt,psfig]{amsart}
\usepackage{amsthm, amsfonts, amssymb,latexsym}
\usepackage{graphics,verbatim}
\usepackage{graphicx}

\newtheorem{mainthm}{Theorem}

\newtheorem{thm}{Theorem}[section]

\newtheorem{lem}[thm]{Lemma}
\newtheorem{prop}[thm]{Proposition}
\theoremstyle{definition}
\newtheorem{defn}[thm]{Definition}
\theoremstyle{remark}

\numberwithin{equation}{section}

\newcommand{\R}{\mathbb R}
\newcommand{\Sph}{\mathbb S}
\newcommand{\C}{\mathbb C}
\newcommand{\N}{\mathbb N}
\newcommand{\CP}{{\mathbb C}{\mathbb P}}
\newcommand{\T}{\mathbb T}
\newcommand{\Z}{\mathbb Z}

\newcommand{\calA}{\mathcal{A}}

\newcommand{\calH}{\mathcal{H}}
\newcommand{\calP}{\mathcal{P}}
\newcommand{\calT}{\mathcal{T}}
\newcommand{\calU}{\mathcal{U}}

\newcommand{\Log}{\mbox{Log}}
\newcommand{\val}{\mbox{val}}

\begin{document}

\title[]{Calabi-Yau Components in general type Hypersurfaces}%
\author{Naichung Conan Leung}%
\address{Department of Mathematics and The Institute of Mathematical Sciences\\
The Chinese University of Hong Kong\\
Shatin, N.T., Hong Kong} \email{leung@ims.cuhk.edu.hk}

\author{Tom Y. H. Wan}%
\address{Department of Mathematics\\
The Chinese University of Hong Kong\\
Shatin, N.T., Hong Kong}%
\email{tomwan@math.cuhk.edu.hk}%
\thanks{Both authors are partially supported by
   Earmarked Grants of Hong Kong CUHK403105}


\begin{abstract}
For a one-parameter family $(V,\{\Omega_i\}_{i=1}^{p_g})$ of general type hypersurfaces with bases of holomorphic $n$-forms, we construct open covers $V=\bigcup_{i=1}^{p_g}U_i$ using tropical geometry. We show that after normalization, each $\Omega_i$ is approximately supported on a unique $U_i$ and such a pair approximates a Calabi-Yau hypersurface together with its holomorphic $n$-form as the parameter becomes large. We also show that the Lagrangian fibers in the fibration constructed by Mikhalkin \cite{Mi03} are asymptotically special Lagrangian. As the holomorphic $n$-form plays an important role in mirror symmetry for Calabi-Yau manifolds, our results is a step toward understanding mirror symmetry for general type manifolds.
\end{abstract}
\maketitle
\setlength{\baselineskip}{20pt}

\section{Introduction}
Calabi-Yau manifolds are K\"ahler manifolds with zero first Chern
class. By Yau's theorem \cite{Yau}, they admit Ricci flat K\"ahler
metrics. They play important roles in String theory as internal
spaces. Up to a scalar multiple, there exists a unique holomorphic
volume form $\Omega\in H^{n,0}\left(  Y\right)  $ on any Calabi-Yau
manifold $Y$. In the SYZ proposal \cite{SYZ} for the Mirror Symmetry
conjecture, Strominger, Yau and Zaslow conjectured that mirror
symmetry is a generalization of the Fourier-Mukai transformation
along dual speacial Lagrangian torus fibrations on mirror Calabi-Yau
manifolds and it is called the ``SYZ transformation". Recall that a
Lagrangian submanifold $L$ in $Y$ is called {\em special} if $\left.
\operatorname{Im}  \Omega \right|_{L}=0$. It is not easy to
construct special Lagrangian fibrations on Calabi-Yau manifolds.
Nevertheless, Lagrangian fibrations do exist on Calabi-Yau
hypersurfaces in $\mathbb{CP}^{n+1}$, or other toric varieties, by
the work of Gross \cite{Gross}, Ruan \cite{Ruan} and others.

There are generalizations of the Mirror Symmetry conjecture for Fano
manifolds (i.e. positive first Chern class) and also recently for
general type manifolds (e.g. negative first Chern class). There are
many Fano manifolds which are toric varieties and therefore they
admit natural Lagrangian torus fibrations. They have canonical
holomorphic volume forms $\Omega$ outside singular fibers which make
the toric fibrations special. The SYZ transformation along these
special Lagrangian fibrations on Fano toric manifolds was studied by
Chan and the first author in \cite{CL}.

This paper is an initial step in our studies of the SYZ mirror
transformation for general type manifolds. In dimension one, for
every $g\ge 2$, there is a family of genus $g$ Riemann surfaces
$V_t$ which degenerate to a connected sum of $g$ copies of elliptic
curves as $t$ goes to infinity, i.e., $V_{\infty}=Y_1\cup\cdots\cup
Y_g$ with each $Y_i$ a smooth elliptic curve. Furthermore, we can
find a base $\Omega_{1,t},\ldots,\Omega_{g,t}$ of $H^{1,0}(V_t)$
such that for each $i\in\{1,\ldots,g\}$, $\Omega_{i,t}$ converges to
a holomorphic volume form on $Y_i$ as $t$ goes to infinity (see
subsection \S\ref{sec-2d}).

In higher dimensions, we cannot expect to have a connected sum
decomposition for general type manifolds $V_t$ into Calabi-Yau
manifolds. Instead, we will show in our main theorem that there is a
basis $\{ \Omega_{1,t}\ldots,\Omega_{p_g,t}\}$ of $H^{n,0}(V_t)$ and
a decomposition
$$
V_t=\bigcup_{i=1}^{p_g}U_{i,t}
$$
such that each $\Omega_{i,t}$ is roughly supported on corresponding
$U_{i,t}$ and $\left( U_{i,t},\Omega_{i,t}\right)$ approximates a
Calabi-Yau manifold $Y_{i,t}$ together with its holomorphic volume
form $\Omega_{Y_{i,t}}$ as $t$ goes to infinity. This is not a
connected sum decomposition as different $U_{i,t}$'s can have large
overlaps. However, it still enables us to have a proper notion of
special Lagrangian fibrations on $V_t$ and study the SYZ
transformation along them.

If $V_{t}$ is a family of general type hypersurfaces in
${\CP}^{n+1}$, i.e. the common degree $d$ of the family of defining
polynomials of $V_{t}$ is bigger than $n+2$, then its geometric
genus
$$
p_{g}\left(  V_{t}\right)  =\dim H^{n,0}\left(  V_{t}\right)  =\binom{d-1}{n+1}\geq2 .
$$
In fact $p_{g}\left(  V_{t}\right)  $ equals to the number of
interior lattice points in $\triangle_{d}$, the standard simplex in
$\mathbb{R}^{n+1}$ spanned by $de_1\ldots de_{n+1}$ and the origin,
where $\{e_\alpha\}_{\alpha=1}^{n+1}$ is the standard basis of
$\R^{n+1}$. That is, if we denote the set of interior lattice points
of $\triangle_{d}$ by $\triangle_{d,\Z}^0$, then
$$
p_{g}\left(  V_{t}\right)  =\# \triangle_{d,\Z}^0.
$$
The analog formula for $p_{g}$ holds true for smooth hypersurfaces
in toric varieties \cite{Fulton}. In this article, we prove the
following
\begin{mainthm}[Main Theorem]
For any positive integers $n$ and $d$ with $d\ge n+2$, there exists
a family of smooth hypersurfaces $V_t\subset \CP^{n+1}$ of degree
$d$ such that $V_t$ can be written as
$$
V_t=\bigcup_{i\in\triangle_{d,\Z}^0} U_{i,t}
$$
where $U_{i,t}$ is a family of open subsets $U_{i,t}\subset V_t$
such that after the normalization $H_t:(\C^*)^{n+1} \to (\C^*)^{n+1}$ defined by
$$ H_t(z_1,\ldots,z_{n+1})=\left( |z_1|^{\frac{1}{\log
t}}\frac{z_1}{|z_1|},\ldots,|z_{n+1}|^{\frac{1}{\log
t}}\frac{z_{n+1}}{|z_{n+1}|} \right),
$$
\begin{enumerate}
\item $U_{i,t}$ is close in Hausdroff distance on
$(\C^*)^{n+1}$ to an open subset of a Calabi-Yau
hypersurface $Y_{i,t}$ in $\CP^{n+1}$.

\item there exists a basis $\{\Omega_{i,t}\}_{i\in\triangle_{d,\Z}^0}$ of $H^{n,0}(V_t)$
such that for each ${i\in\triangle_{d,\Z}^0}$, $\Omega_{i,t}$ is
non-vanishing and close to the holomorphic volume form
$\Omega_{Y_{i,t}}$ of $Y_{i,t}$ on $U_{i,t}$ with respect to the
pull-back metric $H^*_t(g_0)$ of the invariant toric metric $g_0$ on
$(\C^*)^{n+1}$;

\item for any compact subset $B\subset (\C^*)^{n+1}\setminus U_{i,t}$, $\Omega_{i,t}$
tends to zero in $V_t\cap B$ uniformly with respect to
$H^*_t(g_0)$.
\end{enumerate}
\end{mainthm}

\begin{figure}[h]\label{Figure 1}
\centerline{\includegraphics[width=10cm]{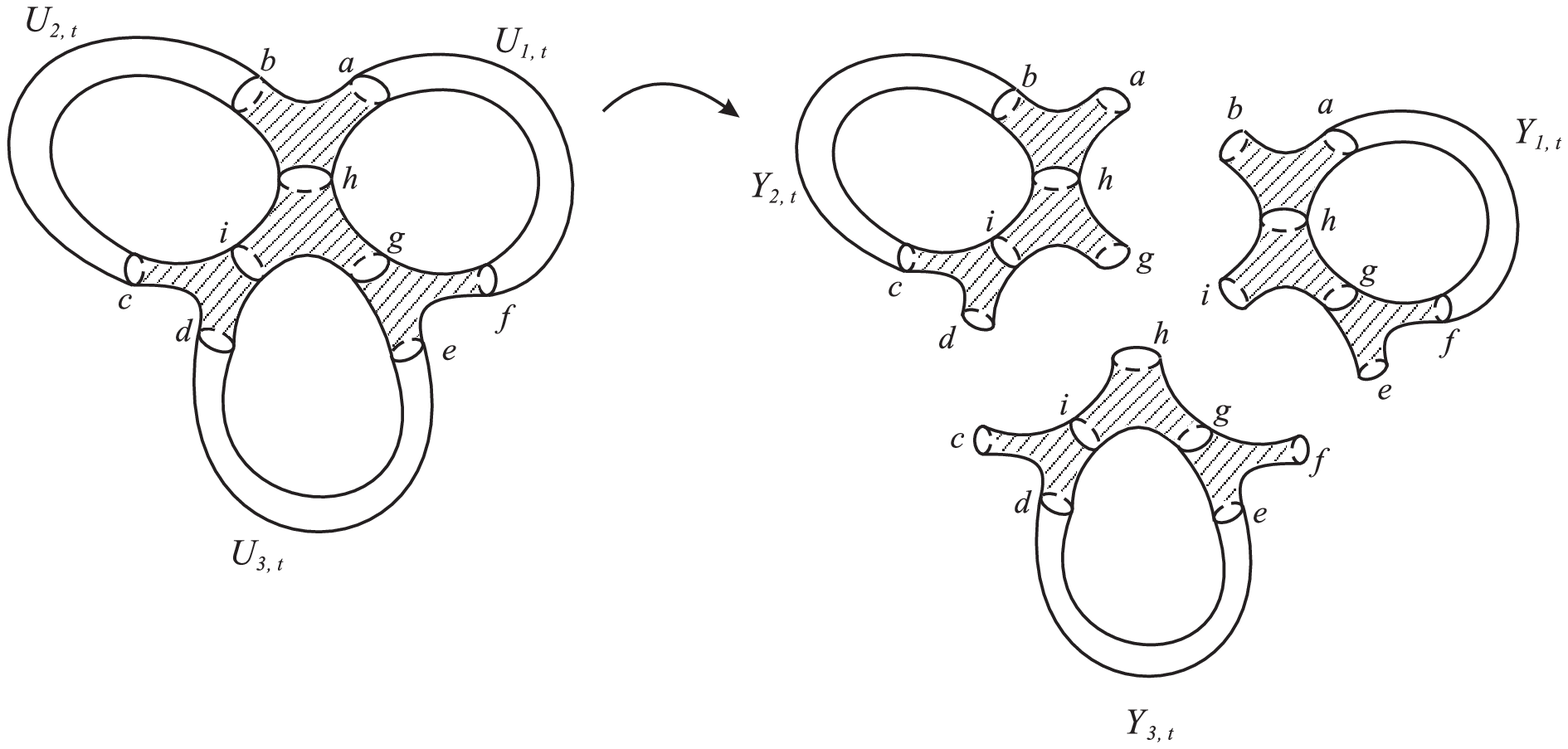}}
\end{figure}
\begin{center}
{\footnotesize Figure 1}
\end{center}

Our proof bases on the results of Mikhalkin \cite{Mi03}. In his
paper, Mikhalkin constructed torus fibrations on general type
hypersurfaces $V$ in $\mathbb{CP}^{n+1}$ and he showed that some of
these fibers are Lagrangian. The technique he employed is tropical
geometry. He constructed a degenerating family $V_{t}$ of
hypersurfaces to decompose $V_{t}$ into union of pairs-of-pants and
also his fibration can be seen from this tropical degeneration. We
are going to make use of his decomposition to construct our open sets
$U_{i,t}$ in the main theorem.

Roughly speaking, the main theorem says that as $t$ approaches
infinity, $V_{t}$ decomposes into $p_{g}$ different Calabi-Yau
manifolds $Y_{i,t}$ and each support a holomorphic $n$-form
$\Omega_{Y_{i,t}}$ on $V_{t}$. (In here, we abused the notion of
``decomposition" since the open sets $U_{i,t}$ that we obtained in
the ``decomposition" do overlap even as $t\to+\infty$.) Therefore we
can speak of \textit{special Lagrangian submanifolds} in $V_{t}$.

\begin{defn} Let $L_{t}\subset V_{t}$ be a smooth family of Lagrangian
submanifolds. We call it  \textit{asymptotically special Lagrangian
} of phase
$\theta$ with respect to the decomposition if for any $\epsilon>0$ we have%
$$
\left|  \left.\operatorname{Im}\left(  e^{\sqrt{-1}\theta}\Omega_{i,t}\right)  \right|_{L_{t}\cap
U_{i,t}}\right|  <\varepsilon
$$
for any $i\in \triangle_{d,\Z}^0$ for sufficiently
large $t$.

If $L_{t}\subset U_{i,t}$ for some $i\in \triangle_{d,\Z}^0$ and $\operatorname{Im}\left(  e^{\sqrt{-1}\theta}%
\Omega_{i,t}\right)  |_{L_{t}}=0$
then we call $L_{t}$ a \textit{special Lagrangian }submanifold in
$V_{t}$.
\end{defn}

In section \S\ref{sec-asl} we show that the Lagrangian fibers in the
torus fibration on $V_{t}$ constructed by Mikhalkin in \cite{Mi03}
are asymptotically special Lagrangians.

We start our proof with some preliminaries on the tropical geometry,
especially on the theorems of Einsiedler-Kapranov-Lind \cite{EKL}
and Mikhalkin \cite{Mi03}. The main results and their proofs will be
stated in section \S\ref{sec-main} and \S\ref{sec-lem}.

\noindent{\sc Acknowledgement} The authors thank the valuable
discussions with Mark Gross and Wei Dong Ruan. We also would like to
thank Miss Suki Chan and Pauline Chan for their helps in making the
figures.

\section{Preliminaries}
\subsection{Amoebas and Viro's patchworking}
Let $V^o$ be a smooth hypersurface in $(\C^*)^{n+1}\subset\CP^{n+1}$
or other toric varieties defined by a Laurent polynomial
$$
f(z)=\sum_j{a_j}z^j,
$$
where $j=(j_1,\ldots,j_{n+1})\in\Z^{n+1}$ are multi-indices. Recall
that the {\em Newton polyhedron} $\triangle\subset\R^{n+1}$ of $f$,
or of $V^o$, is the convex hull in $\R^{n+1}$ of all $j\in\Z^{n+1}$
such that $a_j\neq 0$. According to \cite{Gelfand}, the {\em amoeba}
of $V^o$ is the image
$$
\Log(V^o)\subset \R^{n+1}
$$
under the map $\Log: (z_1,\ldots,z_{n+1})\mapsto
(\log|z_1|,\ldots,\log|z_{n+1}|)$.

In this paper, we are looking for deformation of complex structures
on $V^o$ together with corresponding basis of holomorphic $n$-forms
satisfying a special limiting property. This leads us to consider
deformation of the polynomial $f$ used by the {\em Viro's
patchworking} \cite{Viro} and non-Archimedean amoeba .

Let $v: \triangle_{\Z} \to \R$, where $\triangle_{\Z}=\triangle \cap\Z^{n+1}$, be any function and
$f(z)=\displaystyle\sum_{j\in\triangle_{\Z}} a_jz^j$, $a_j\neq 0$ for any
$j\in\triangle_{\Z}$, be any polynomial. The {\em
patchworking polynomial} is defined for all $t>0$ by
$$
f_t^v(z)=\sum_{j\in\triangle_{\Z}} a_jt^{-v(j)}z^j.
$$
The family $f_t^v$ can be treated as a single polynomial in
$(K^*)^{n+1}$, where $K^*=K\setminus\{0\}$ and $K$ is the field of
Puiseux series with complex coefficients in $t$. In order to match
the notation in the literatures, for instance \cite{EKL}, we set
$\tau=t^{-1}$ and let
$$\C((\tau^{q}))=\left\{g(\tau^{q})=\sum_{k=m}^\infty g_m (\tau^q)^k \right\}$$
be the field of formal (semi-finite) Laurent series in $\tau^q$.
Then the field of Puiseux series is
$$
K=\bigcup_{m\ge 1}\C((\tau^{\frac{1}{m}})).
$$
The field $K$ is algebraically closed \cite{Cohn} and has a
valuation defined by
$$
\val_K\left( \sum_{q\in \Lambda_b}b_q\tau^q \right)=\min\Lambda_b.
$$
for $b=\displaystyle\sum_{q\in \Lambda_b}b_q\tau^q \in K$. It is then easy to see
that the field $K$ can also be represented by the field of Puiseux
series
$$
\tilde{b}=\sum_{p\in \tilde{\Lambda_{{b}}}}\tilde{b}_pt^p
$$
with $\max\tilde{\Lambda}_p<+\infty$ and valuation
$\val_K(\tilde{b})=-\max\tilde{\Lambda}_{{b}}$.

Since $e^{-\val_K}$ defines a norm $\|\cdot\|_K$ on $K$, we can
define $\Log_K$ on $(K^*)^{n+1}$ analog to $\Log$ on $(\C^*)^{n+1}$
by
\begin{eqnarray*}
\Log_K(a_1,\ldots,a_{n+1})&=&(\log\| a_1\|_K,\ldots, \log\|
a_{n+1}\|_K)\\
&=& -(\val_K(a_1),\ldots,\val_K(a_{n+1})).
\end{eqnarray*}
Then for $V_K\subset(K^*)^{n+1}$, the image set
$\calA_K=\Log_K(V_K)$ is called accordingly the {\em (non-Archimedean)
amoeba} of $V_K$. It is clear that $\calA_K=-\calT(V_K)$, where
$\calT(V_K)$ is the tropical variety of $V_K$ which is defined as
the closure of $\val_K(V_K)$ \cite{EKL}.

Note that for our family $f_t^v(z)=\displaystyle\sum_{j\in\triangle_{\Z}}
a_jt^{-v(j)}z^j$, the coefficient of $z^j$ is $a_jt^{-v(j)}\in K$
which has valuation $\val_K(a_jt^{-v(j)})=v(j)$. This match the
convention in \cite{Mi03}.

Following the construction of \cite{Mi03}, for a finite
set $A$ in $\Z^{n+1}$ and a real valued function $v:A\to\R$ on $A$, one defines
 $\Pi_v$ to be the set of non-smooth points (called {\em corner locus} in
\cite{Mi03}) of the Legendre transform $L_v:\R^{n+1}\to \R$ of $v$. Here $L_v(x)$ is
defined by
$$
L_v(x)=\max_{i\in A}l_{v,i}(x),
$$
where $l_{v,i}(x)=\langle x, i\rangle -v(i)$ with $\langle\cdot\, ,
\cdot\rangle$ is the standard inner product on $\R^{n+1}$. In
particular, the interior of a top dimensional face of $\Pi_v$ is
given by
$$
\frak{F}(j^{(1)},j^{(2)})=\{
x\in\R^{n+1}\,:\,l_{v,j^{(1)}}(x)=l_{v,j^{(2)}}(x)>l_{v,j}(x),\,
\forall\, j\neq j^{(1)},\, j^{(2)}\}.
$$
It was proved in
\cite{Mi03} that $\Pi_v$ is a balanced polyhedral complex dual to
certain lattice subdivision of the convex hull $\triangle$ of $A$ in
$\R^{n+1}$. We refer the reader to \cite{Mi03} or the appendix for
the definition of a balanced polyhedral complex. We have the
following result of Einsiedler-Kapranov-Lind \cite{EKL}.
\begin{thm}
If $V_K\subset (K^*)^{n+1}$ is a hypersurface given by a polynomial
$f=\displaystyle\sum_{j\in\triangle_\Z} a_jz^j$, $a_j\in K^*$. Then
the (non-Archimedean) amoeba $\calA_K$ of $V_K$ is the balanced
polyhedral complex $\Pi_v$ corresponding to the function
$v(j)=\val_K(a_j)$ defined on the lattice points of the Newton
polyhedron $\triangle$.
\end{thm}
Note that the theorem in \cite{EKL} is originally stated for the
tropical variety $\calT(V_K)$ instead of $\calA_K$.

Now we can describe the limiting behavior of the family of varieties
$V^o_t=\{f^v_t=0\}$ in $(\C^*)^{n+1}$ as $t\to+\infty$. For each
$t>0$, we define the amoeba of $V^o_t$ with respect to $t$ by
$$
{\calA}_t=\Log_t(V_t^o)\subset \R^{n+1},
$$
where
$\Log_t(z_1,\ldots,z_{n+1})=(\log_t|z_1|,\ldots\log_t|z_{n+1}|)$ on
$(\C^*)^{n+1}$, where $\log_t\tau=\log\tau/\log t$ for $\tau>0$. If
we denote accordingly $\calA_K=\Log_K(V_K)$ the non-Archimedean
amoeba of the family $f_t^v$ regarded as a single polynomial in the
field $K$ of Puiseux series. Then, we have the following theorem of
Mikhalkin \cite{Mi03} which is needed in the proofs of our
assertions.
\begin{thm}
The amoebas $\calA_t$ converge in the Hausdorff distance on
$\R^{n+1}$ to the non-Archimedean amoeba $\calA_K$ as $t\to+\infty$.
\end{thm}
Recall that the Hausdorff distance between two closed subsets $A$
and $B$ in $\R^{n+1}$ is given by
$$
d_{H}(A,B)=\max\left\{ \sup_{a\in A}d_{\R^{n+1}}(a,B),\,\sup_{b\in
B}d_{\R^{n+1}}(A,b)\right\}.
$$

\subsection{Maximal dual complex}
As we mentioned, it was proved in \cite{Mi03} that $\Pi_v$ is a
balanced polyhedral complex dual to certain lattice subdivision of
the convex hull $\triangle$ of $A$ in $\R^{n+1}$. In general, any
n-dimensional balanced polyhedral complex $\Pi$ in $\R^{n+1}$
determines a convex lattice polyhedron $\triangle\subset\R^{n+1}$
and a lattice subdivision of $\triangle$. We call $\Pi$ a
{\em maximal} polyhedral complex if the elements of the subdivision
are simplices of volume $\frac{1}{(n+1)!}$, i.e, the corresponding
subdivision is a unimodular lattice triangulation. Note that not all
convex lattice polyhedron admit unimodular lattice triangulation.
Therefore, not all convex lattice polyhedron admit maximal dual
complex. If it does, then we have the following result of
\cite{Mi03}.
\begin{prop}\label{prop-topPi}
If $\Pi$ is a maximal dual $\triangle$-complex, then $\Pi$ is
homotopy equivalent to the bouquet of $\#\triangle_{\Z}^0$ copies of
$\Sph^n$, where
$\triangle_{\Z}^0=\left(\operatorname{Int}\triangle\right)\cap\Z^{n+1}$
is the set of interior lattice points of $\triangle$.
\end{prop}
However, the converse of the proposition is not true. A non-maximal
dual $\triangle$-complex may still have the homotopy type stated in
the proposition.

It was also shown in \cite{Mi03} that on each maximal complex
$\Pi\subset \R^{n+1}$, there is a canonical choice of cutting locus
$\Xi$ such that each connected component $\calU_k$, ($k=1,\ldots, l$) called {\em
primitive piece}, of $\Pi\subset \Xi$ is equivalent to an open neighborhood of the vertex in the {\em primitive complex} $\Sigma_n\subset \R^{n+1}$ which is the set of non-smooth points of
the function $H(x_1,\ldots,x_{n+1})=\max\{0,x_1,\ldots,x_{n+1}\}$.
That is, there exists $M_k\in ASL_{n+1}(\Z)=SL_{n+1}(\Z) \ltimes
\Z^n$ such that $M_k(\calU_k)$ is an open set of $\Sigma_n$
containing the vertex. Furthermore, these open sets are parametrized
by the vertices of $\Pi$. Since $\Pi$ is dual to the lattice subdivision of $\triangle$ with simplices of volume $\frac{1}{(n+1)!}$, we must have exactly $(n+1)!\mbox{vol}(\triangle)$ distinct $\calU_k$, i.e. $l=(n+1)!\mbox{vol}(\triangle)$
\pagebreak
\begin{figure}[h]\label{Figure 2}
\centerline{\includegraphics[width=5cm]{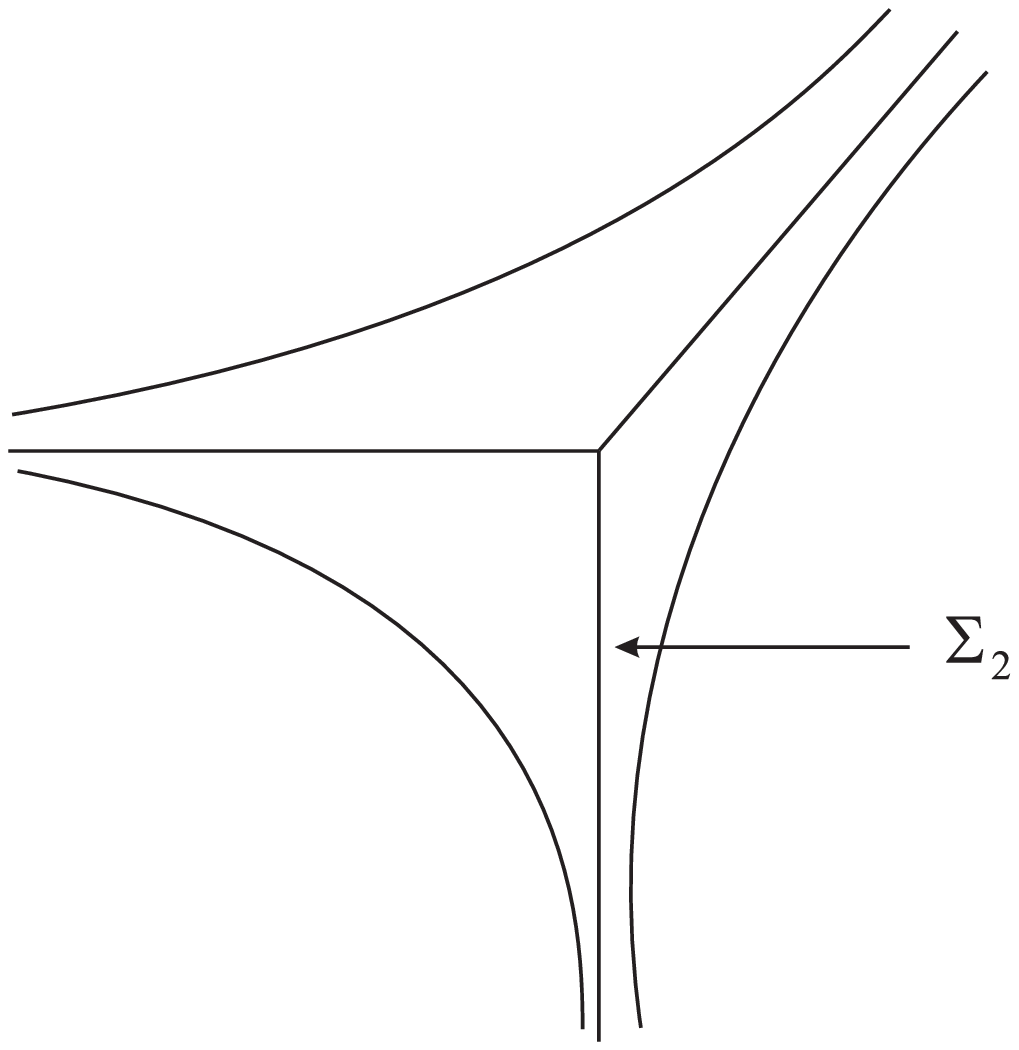}}
\end{figure}
\begin{center}
{\footnotesize Figure 2: The primitive complex $\Sigma_2$ and the amoeba of $1+z_1+z_2=0$.}
\end{center}


\subsection{Pairs-of-pants decomposition and stratified fibration}
In this subsection, we state the pairs-of-pants decomposition and
existence of stratified fibration theorem of Mikhalkin \cite{Mi03}
which is the main ingredient of the proof of our results. We start
with the definition of pair-of-pants and stratified fibration given
in \cite{Mi03}.

As in \cite{Mi03}, we denote by $\calH$ a union of $n+2$ generic
hyperplanes in $\CP^n$ and $\calU$ the union of the corresponding
$\epsilon$-neighborhoods for a small $\epsilon>0$. Then
$\overline{\calP_n}=\CP^n\setminus \calU$ is called the {\em
$n$-dimensional pair-of-pants} while ${\calP_n}=\CP^n\setminus
\calH$ the {\em $n$-dimensional open pair-of-pants}. It is clear
that $\calP_1$ is diffeomorphic to a 2-sphere with three punctures
and $\overline{\calP_1}$ is diffeomorphic to a 2-sphere with three
holes, or equivalently, a closed disk with two holes. That is, the
definition is a generalization of the classical pair-of-pants in one
complex dimension.

If $V$ and $F$ are smooth manifolds, and $\Pi$ is a maximal dual
$\triangle$-complex of a lattice polyhedron $\triangle$ of full
dimension in $\R^{n+1}$, then a smooth map $\lambda:\,V\to\Pi$ is a
{\em stratified $F$-fibration} if it satisfies
\begin{enumerate}
\item the restriction of $\lambda$ over each open $n$-cell
$e\subset\Pi$ is a trivial fibration with fiber $F$;
\item for each pair of integers $(l,k)$ with $0\le k \le l \le n$,
there exists a smooth ``model" map depending only on $l$ and $k$,
$\lambda_{l,k}:\,V_{l,k}\to\Pi_{l,k}$ with $\Pi_{l,k}$ diffeomorphic
to $\R^k\times\Sigma_{l-k}\times [0,+\infty)^{n-l}$ such that any
$(l,k)$-point of $\Pi$ has a neighborhood $U$ such that
$$
\lambda|_U :\,\lambda^{-1}(U)\to U
$$
is diffeomorphic to the model map.
\end{enumerate}

Now, we can state the pairs-of-pants decomposition and existence of
stratified fibration theorem of Mikhalkin \cite{Mi03}.
\begin{thm}\label{thm-pop}
Let $V$ be a smooth hypersurface in $\CP^{n+1}$ defined by a
polynomial with Newton polyhedron
$\triangle_d$.
Then for every maximal dual $\triangle_d$-complex $\Pi$, there
exists a stratified $\T^n$-fibration $\lambda:\, V\to \Pi$
satisfying
\begin{enumerate}
\item the induced map $\lambda^*: H^n(\Pi, \Z)\approx \Z^{p_g} \to H^n(V,\Z)$ is
injective, where $p_g=h^{n,0}(V)$ is the geometric genus of $V$;
\item for each primitive piece $\calU_k$ of $\Pi$, $\lambda^{-1}(\calU_k)$
is diffeomeorphic to an open pair-of-pants $\calP_n$;
\item for each $n$-cell $e$ of $\Pi$, there exists a point $x\in e$
such that the fiber $\lambda^{-1}(x)$ is a Lagrangian $n$-torus
$\T^n\subset V$;
\item there exists Lagrangian embedding $\phi_i:\Sph^n\to V$,
$i=1,\ldots, p_g$ such that the cycles $\lambda\circ\phi_i(\Sph^n)$
form a basis of $H_n(\Pi)$.
\end{enumerate}
\end{thm}

\pagebreak
\begin{figure}[h]\label{Figure 3}
\centerline{\includegraphics[width=10cm]{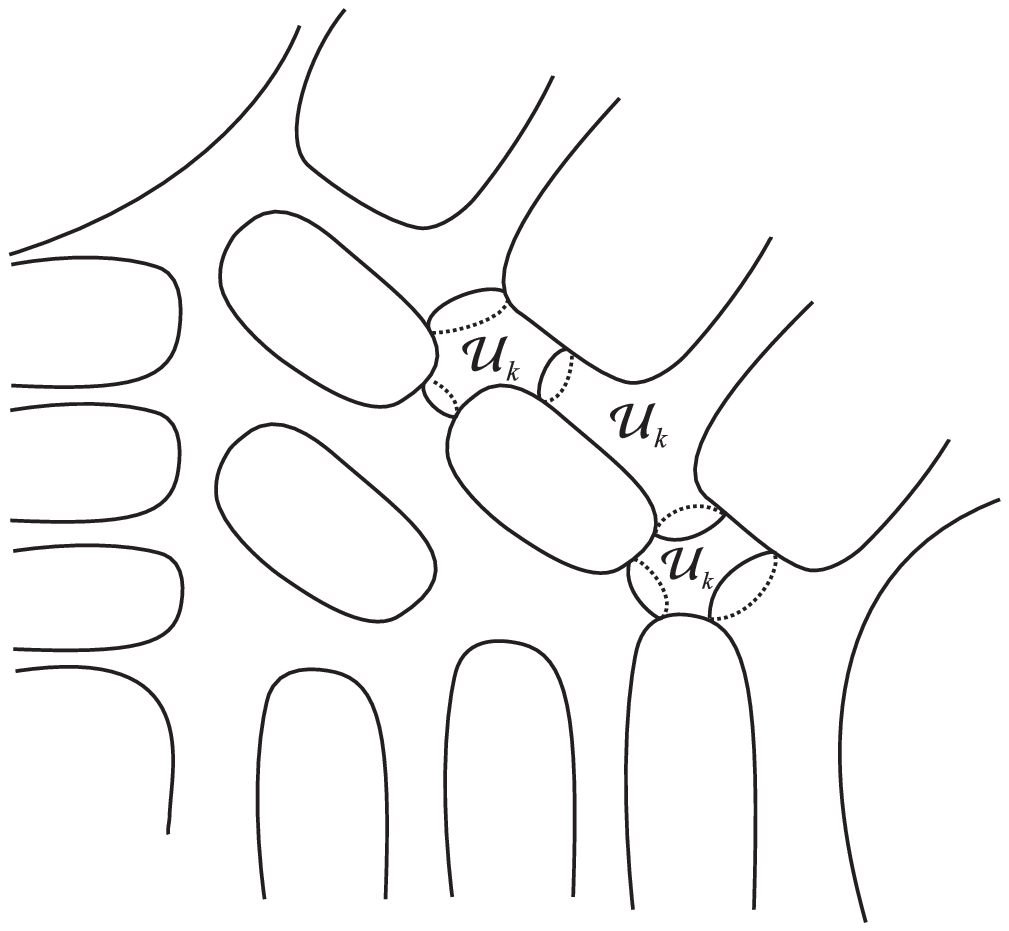}}
\end{figure}
\begin{center}
{\footnotesize Figure 3: Illustration of some pairs-of-pants $\calU_k$.}
\end{center}

\subsection{Key lemma}
To prove the main theorem, we need to show the existence of a real
valued function $v:\triangle_{d,\Z}\to \R$, where
$\triangle_{d,\Z}=\triangle_{d}\cap\Z^{n+1}$ and $d\ge n+2$, such
that the corresponding lattice subdivision of $\triangle_d$ dual to
the balanced polyhedral complex $\Pi_v$ satisfies some special properties
which are needed to obtain the ``decomposition" into Calabi-Yau pieces claimed in our main
result. The special properties that we need are the property (2) in the
following lemma. Existence of function $v$ with property (1) only is
well-known.
\begin{lem}\label{lem-v}
Let $\triangle_{d}$, $d\ge n+2$, be the simplex in $\R^{n+1}$ with
vertices $\{0,de_1,\ldots, de_{n+1}\}$, where
$\{e_1,\ldots,e_{n+1}\}$ is the standard basis of $\R^{n+1}$, and
the function $v:\triangle_{d,\Z} \to \R$ be defined by
$$
v(j_1,\ldots,j_{n+1})=\sum_{\alpha=1}^{n+1} j_\alpha^2+\left(
\sum_{\alpha =1}^{n+1} j_\alpha \right)^2.
$$
Then
\begin{enumerate}
\item the balanced polyhedral complex $\Pi_v$ corresponding to $v$ is
a maximal dual complex of $\triangle_{d}$;
\item the corresponding
subdivision of $\triangle_{d}$ has the property that for each
point $i=(i_1,\ldots, i_{n+1})\in \triangle_{d,\Z}^0$, the lattice subdivision restricts to
a lattice subdivision of the translated simplex
$i-{\i}+\triangle_{n+2}$, where ${\i}=(1,\ldots,1)\in\Z^{n+1}$.
\end{enumerate}
\end{lem}

\begin{figure}[h]\label{Figure 4}
\centerline{\includegraphics[width=8cm]{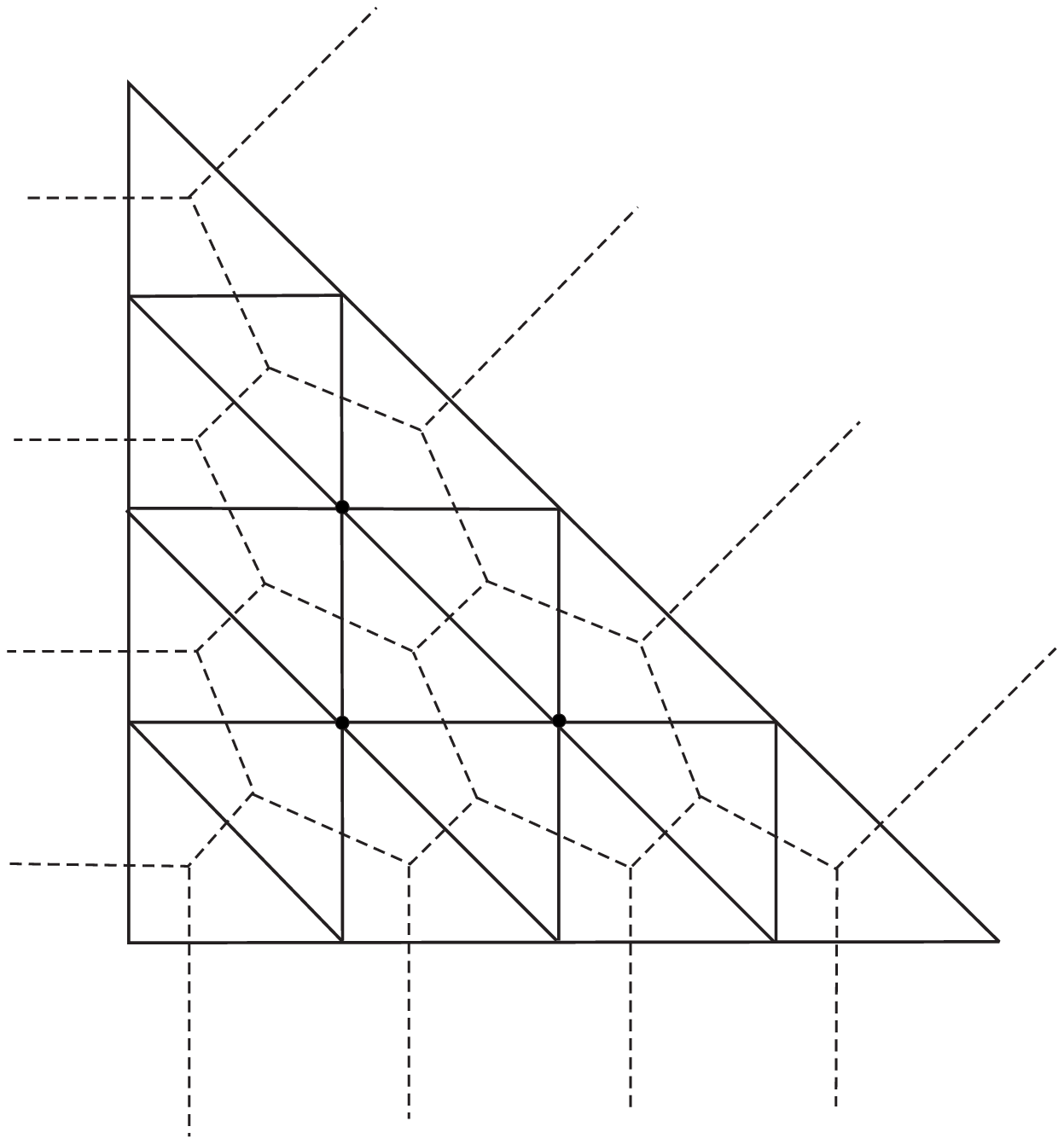}}
\end{figure}
\begin{center}
{\footnotesize Figure 4: A maximal dual $\triangle_4$ complex with the required property.}
\end{center}

Note that we do not claim that the subdivision is translational
invariant which probably is not true. We only claim that for each
interior lattice point, one can find a union of simplexes which form
a standard simplex with exactly one interior lattice point. This
clearly corresponds to the Calabi-Yau situation. The proof of this
lemma will be given in the last section of this article. 


\section{Proof of the main theorem}\label{sec-main}
In this section, we give the proof of the main theorem.  Recall that
we are free to use any hypersurface defined by a homogeneous
polynomial of degree $d$ of $n+1$ variables to replace $V$ in order
to describe $V$ as a smooth manifold or as a symplectic manifold.
The idea of tropical geometry leads us to consider the submanifold
$V^o=V\cap{(\C^*)}^{n+1}$. Since all smooth hypersurfaces with the same Newton polyhedron are
isotopic, we are free to choose the non-zero coefficients. In
particular, we may take $a_j=1$ for all $j\in\triangle_{\Z}$.
Therefore, a deformation of complex
structures on $V$ can be given by a patchworking polynomial
$f_t(z)=\displaystyle\sum_{j\in\triangle_\Z}t^{-v(j)}z^j$, $t>0$ in $z\in{(\C^*)}^{n+1}$. We assume
$f_t(z)$ is generic with Newton polyhedron
$\triangle=\triangle_d$
and $v(j)<\infty$ for all $j\in \triangle_{d,\Z}$. That is,
$f_t(z)$ contains all possible monomials of degree less than or
equal to $d$.


The result of Einsiedler-Kapranov-Lind in the previous section
states that the amoebas $\calA_t=\Log_t(V_t^o)$ converge in the
Hausdorff distance on $\R^{n+1}$ to the non-Archimedean amoeba
$\calA_K=\Pi_v$ as $t\to+\infty$. By construction, the top
dimensional faces are given by two maximal terms in $f_t(z)$ as
$t\to +\infty$. To be precise, we note that the highest exponent of
$t$ for the term $t^{-v(j)}z^j$ in the polynomial $f_t$ is given by
$$
l_{v,j}(x)=\langle j, x\rangle - v(j),
$$
where $x=\displaystyle\lim_{t\to+\infty}\Log_t(z)=
\displaystyle\lim_{t\to+\infty}(\log_t|z_1|, \ldots,
\log_t|z_{n+1}|) \in \R^{n+1}$. Then, the interior of a top
dimensional face of $\calA_K=\Pi_v$ associated to the terms
$t^{-v(j^{(1)})}z^{j^{(1)}}$ and $t^{-v(j^{(2)})}z^{j^{(2)}}$ is
$$
\frak{F}(j^{(1)},j^{(2)})=\{
x\in\R^{n+1}\,:\,l_{v,j^{(1)}}(x)=l_{v,j^{(2)}}(x)>l_{v,j}(x),\,
\forall\, j\neq j^{(1)},\, j^{(2)}\}.
$$

By the results of tropical geometry \cite{Mi03}, the generators of the
$n$-dimensional homology of the amoeba are exactly given by the limit
of the boundaries of the domains on which a term corresponds to an
interior lattice point of the Newton polyhedron is maximal. Recall that the set of interior lattice
points of $\triangle_d$ is exactly equal to $p_g=\left(
\begin{array}{c} d-1 \\n+1 \end{array} \right)$ (see for instance \cite{Fulton}) and note that
$\frak{F}(i,j)\neq \emptyset$ only when $i$ and $j$ is connected by
an edge in the lattice subdivision of $\triangle_d$ dual to $\Pi_v$.
Then for each $i\in\triangle_{d,\Z}^0$,
$C^{\wedge}_i=\bigcup_j\overline{\frak{F}(i,j)}$ forms an $n$-cycle
representing an element of $H_n(\Pi_v,\Z)$. It is clear from
Mikhalkin's results mentioned in the previous section that the
classes $\{[C^{\wedge}_i]\}_{i\in\triangle_{d,\Z}^0}$ are the generators of
$H_n(\Pi_v,\Z)$. Since $\Pi_v=\calA_K$, this gives the generators of
$H_n(\calA_K,\Z)$ and hence $H_n(\calA_t,\Z)$ for large $t$. From
this observation, we first prove the following result which gives
partial results of the main theorem.

\begin{thm}\label{thm-form}
Let $v:\triangle_{d,\Z}\to \R$ be a real valued function such that
the set $\Pi_v$ of non-smooth points of the Legendre transform of
$v$ is a maximal dual $\triangle_d$-complex, $d\ge n+2$. Let
$f_t=\displaystyle\sum_{j\in\triangle_{d,\Z}}t^{-v(j)}z^j$ ($t>0$)
be the patchworking polynomial of degree $d$ defined by $v$ with
non-Archimedean amoeba $\calA_K=\Pi_v$. Denote
$V_t=\{f_t=0\}\subset\CP^{n+1}$. Then for all $t>0$, there exists a
basis $\{\Omega_{i,t}\}_{i\in\triangle_{d,\Z}^0}$ of $H^{n,0}(V_t)$,
and open subsets ${U}^{\wedge}_{i,t} \subset V_t\cap (\C^*)^{n+1}$
such that for each $i\in\triangle_{d,\Z}^0$,

\begin{enumerate}
\item $\Log_t({U}^{\wedge}_{i,t})$ tends to an $n$-cycle $C^{\wedge}_i$ such that
$\{[C^{\wedge}_i]\}_{i\in\triangle_{d,\Z}^0}$ forms a basis of
$H_n(\calA_K(V))$,

\item $\Omega_{i,t}$ is nonvanishing on ${U}^{\wedge}_{i,t}$ for large $t$, and

\item for any compact subset $B\subset \CP^{n+1}\setminus {U}^{\wedge}_{i,t}$, $\Omega_{i,t}$
tends to zero in $V_t^o\cap B$ uniformly with respect to the metric
induced from the pull-back metric $H^*_t(g_0)$ of the invariant
metric of the torus $(\C^*)^{n+1}$.
\end{enumerate}
\end{thm}

\pagebreak
\begin{figure}[h]\label{Figure 5}
\centerline{\includegraphics[width=10cm]{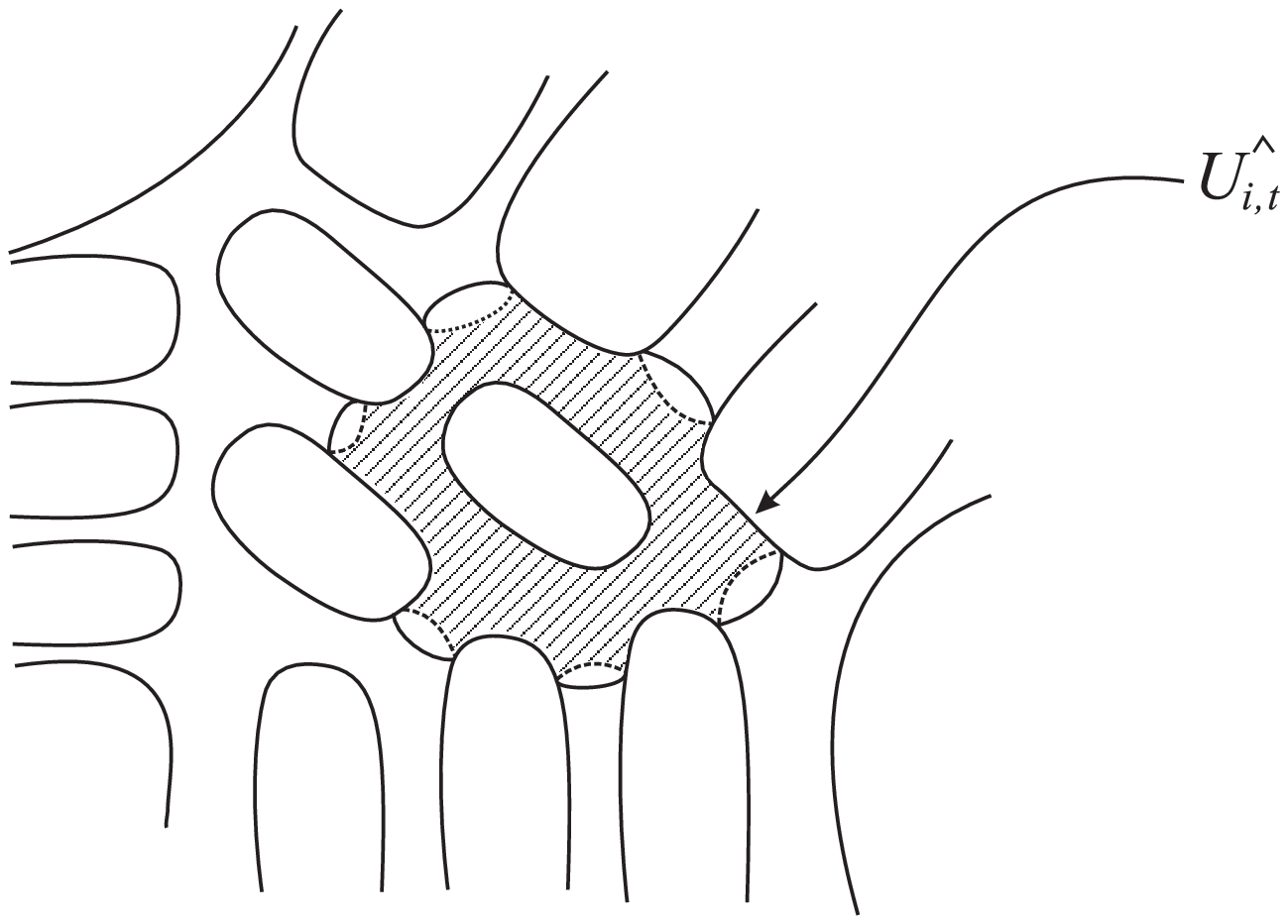}}
\end{figure}
\begin{center}
{\footnotesize Figure 5: Illustration of an open set
$U^{\wedge}_{i,t}$.}
\end{center}

\begin{proof}
By the above observation, we consider for each
$i\in\triangle_{d,\Z}^0$, the cycle
$C^{\wedge}_i=\bigcup_j\overline{\frak{F}(i,j)}$.  Denote the
$1/t$-neighborhood of $C^{\wedge}_i$ in $\Pi_v$ (not in $\R^{n+1}$)
by $C^{\wedge}_{i,t}$. (In fact, one can take any function of $t$
instead, as long as it tends to $0$ as $t\to+\infty$.) Choose
$\delta>0$ in such a way that $\Pi_v$ is a deformation retract of
the $\delta$-tubular neighborhood $T_{\delta}$ of $\Pi_v$. Then for
any $t>0$ such that $\Log_tV_t\subset T_{\delta}$, we define the
open set
$$
{U}^{\wedge}_{i,t}=\lambda^{-1}\left(C^{\wedge}_{i,t}\right),
$$
where $\lambda:V_t\to\Pi_v$ is the stratified $\T^n$-fibration given
by Mikhalkin \cite{Mi03}. Then it is clear that condition (1) is
satisfied, namely $\Log_t({U}^{\wedge}_{i,t})$ tends to
$C^{\wedge}_i$ in Hausdroff distance as $t\to \infty$. This proves
the first statement.

\pagebreak
\begin{figure}[h]\label{Figure 6}
\centerline{\includegraphics[width=6cm]{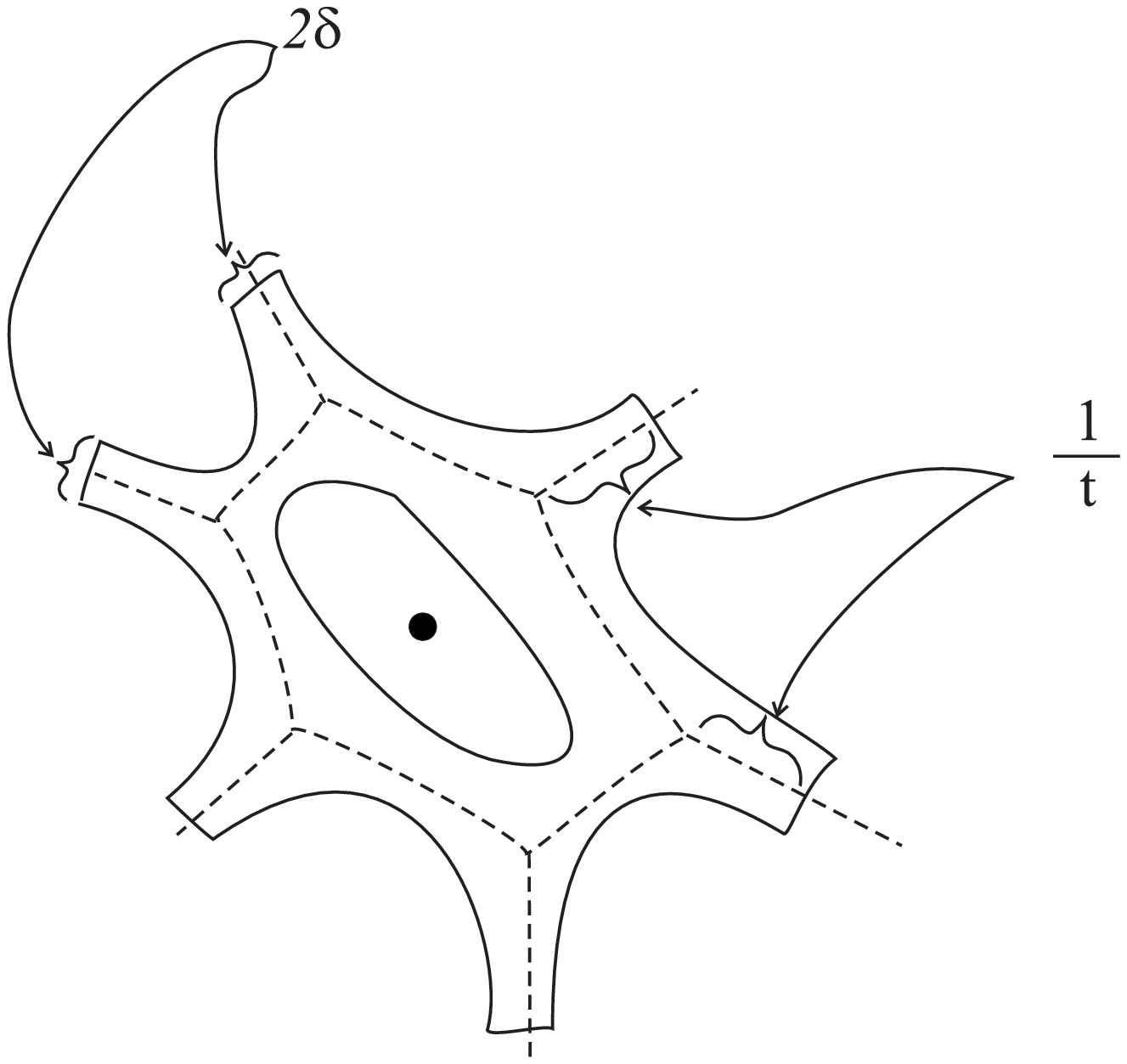}}
\end{figure}
\begin{center}
{\footnotesize Figure 6: The construction of the neighborhood
${C}^{\wedge}_{i,t}$.}
\end{center}

To see the other two statements, we use the well-known fact \cite{Grif} that on
the variety $V_{t}$, the Poincar\'e residues of
$f_t^{-1}dz_1\wedge\cdots\wedge dz_{n+1}$ define a holomorphic
$n$-form on $V_t$; and all elements in $H^{n,0}(V_t)$ are of the
form
$$
\Omega=P(z)\operatorname{Res}\frac{dz_1\wedge\cdots\wedge dz_{n+1}}{f_t(z)},
$$
where $P(z)$ is a polynomial of degree at most $d-(n+2)$. To
simplify notation, we omit the subscribe $t$ in the following
calculations.

For any $\alpha\in\{1,\ldots,n+1\}$, in the region where
$f_{z_{\alpha}}=\frac{\partial f}{\partial z_\alpha}\neq 0$, the residue $\Omega_o$ of
$f(z)^{-1}{dz_1\wedge\cdots\wedge dz_{n+1}}$ is given by
$$
\Omega_o=(-1)^{\alpha-1}\frac{dz_1\wedge\cdots\widehat{dz_\alpha}\cdots\wedge
dz_{n+1}}{f_{z_\alpha}}.
$$
Now, for each interior lattice point $i\in \triangle_{d,\Z}^{0}$, we
define
$$
\Omega_i = (\log t)^{-n}\frac{t^{-v(i)}z^{i}}{z_1\cdots z_{n+1}}
\Omega_o,
$$
Note that $i$ belongs to $\triangle_{d,\Z}^0$ implies $i_\alpha\ge
1$ for all $\alpha=1,\ldots, n+1$. This shows that
${t^{-v(i)}z^{i}}/{z_1\cdots z_{n+1}}$ is a polynomial of degree
less than or equal to $d-1-(n+1)=d-(n+2)$ and hence $\Omega_i$ is a
holomorphic $n$-form on $V_t$. It is also clear from the
construction that $\Omega_i$, $i\in\triangle_{d,\Z}^0$, form a basis
of $H^{n,0}(V_t)$. Explicitly, in the region with $f_{z_\alpha}\neq
0$,
\begin{eqnarray*}
\Omega_i & = &
(-1)^{\alpha-1}(\log t)^{-n}\frac{t^{-v(i)}z^{i}}{z_1\cdots
z_{n+1}}\frac{dz_1\wedge\cdots\widehat{dz_\alpha}\cdots\wedge
dz_{n+1}}{f_{z_\alpha}} \\
& = & (-1)^{\alpha-1}(\log t)^{-n}\frac{t^{-v(i)}z^{i}}{\sum_j
j_\alpha t^{-v(j)}z^j}
\frac{dz_1}{z_1}\wedge\cdots\widehat{\left(\frac{dz_\alpha}{z_\alpha}\right)}
\cdots\wedge\frac{dz_{n+1}}{z_{n+1}}.
\end{eqnarray*}

Now each (non-empty) face $\frak{F}(i,j^{(1)})$ of the $n$-cycle
$C^{\wedge}_i$ of $\calA_K$ corresponding to $i$ is given by
$$
l_{v,i}(x)=l_{v,j^{(1)}}(x)
> l_{v,j}(x),\quad j\neq i,\, j^{(1)},
$$
where ${j^{(1)}}\in(\triangle_{d,\Z}\setminus\{i\})$. Namely,
$t^{-v({i})}z^{{i}}$  and $t^{-v(j^{(1)})}z^{j^{(1)}}$ are the two
maximal terms of $f_t(z)$ determining the face as $t\to+\infty$.
Therefore for any compact subset
$R\Subset\operatorname{Int}(\frak{F}(i,j^{(1)}))$ the terms
$t^{-v({i})}z^{{i}}$ and $t^{-v(j^{(1)})}z^{j^{(1)}}$ dominate other
terms of $f_t$ in a neighborhood of $\lambda^{-1}(R)\subset
V_t\cap(\C^*)^{n+1}$ in $(\C^*)^{n+1}$ as $t\to+\infty$.

For each $\alpha\in \{1,\ldots,n+1\}$, the definition of $f_t(z)$
gives (omitting the subscribe $t$)
\begin{eqnarray*}
z_{\alpha}f_{z_\alpha}&=&\sum_j j_\alpha t^{-v(j)}z^j\\
&=&{i}_\alpha t^{-v({i})}z^{{i}}+ j^{(1)}_{\alpha} t^{-v(j^{(1)})}z^{j^{(1)}} + \cdots \\
&=&{i}_\alpha\left(t^{-v({i})}z^{{i}} + t^{-v(j^{(1)})}z^{j^{(1)}} \right) +
\left( j^{(1)}_{\alpha} - {i}_\alpha \right)
t^{-v(j^{(1)})}z^{j^{(1)}}+
 \cdots\\
&=&{i}_\alpha\left(f+\cdots\right) +
\left( j^{(1)}_{\alpha} - {i}_\alpha \right)
t^{-v(j^{(1)})}z^{j^{(1)}}+
 \cdots
\end{eqnarray*}
where ``$\cdots$" denotes the terms in $f$, up to multiple of a
constant, other than $t^{-v(i)}z^{i}$ and
$t^{-v(j^{(1)})}z^{j^{(1)}}$, and $i_\alpha$ and $j^{(1)}_{\alpha}$ are the
$\alpha$-components of $i$ and $j^{(1)}$ respectively.

Since $j^{(1)}\neq {i}$, there is an index $\alpha\in
\{1,\ldots,n+1\}$ such that $j^{(1)}_{\alpha}\neq {i}_\alpha$.
Therefore, for this $\alpha$ and sufficiently large $t$,
$f_{z_\alpha}\neq 0$ in a neighborhood of $\lambda^{-1}(R)$ for any
compact subset $R\Subset\operatorname{Int}(\frak{F}(i,j^{(1)}))$ of the
interior of the face $\frak{F}(i,j^{(1)})$. Putting the above expression
into the definition of $\Omega_i$ and using $f=0$ on $V_t$, we have
$$
\Omega_i = (\log t)^{-n} \left[
\frac{(-1)^{\alpha-1}t^{-v(i)}z^{i}}{(j^{(1)}_{\alpha}-{i}_\alpha)
t^{-v(j^{(1)})}z^{j^{(1)}}+ \cdots} \right]
\frac{dz_1}{z_1}\wedge\cdots\widehat{\left(\frac{dz_\alpha}{z_\alpha}\right)}
\cdots\wedge\frac{dz_{n+1}}{z_{n+1}}.
$$
Let $\xi=H_t(z)=\left( |z_1|^{\frac{1}{\log
t}}\frac{z_1}{|z_1|},\ldots,|z_{n+1}|^{\frac{1}{\log
t}}\frac{z_{n+1}}{|z_{n+1}|} \right)$ be the normalization mentioned in (2)
of the main theorem. Then $\Log |\xi|=\Log_t |z|$ and
$$
(\log t)^{-1}\frac{dz_\beta}{z_\beta} = \frac{d|\xi_\beta|}{|\xi_\beta|}
+ \sqrt{-1}(\log t)^{-1}d\arg{z_\beta}.
$$
Therefore, by using $t^{-v(j^{(1)})}z^{j^{(1)}} + t^{-v({i})}z^{{i}}+\cdots=f(z)=0$, as $t\to+\infty$
$$
\Omega_i\to \left[
\frac{(-1)^{\alpha-1}}{i_{\alpha}-j^{(1)}_\alpha} \right]
\frac{d|\xi_1|}{|\xi_1|}\wedge\cdots\wedge\widehat{\left(\frac{d|\xi_\alpha|}{|\xi_\alpha|}\right)}
\wedge\cdots\wedge\frac{d|\xi_{n+1}|}{|\xi_{n+1}|}
$$
locally, with respect to the metric induced from the pull-back
$H^{*}_t(g_0)$ of the invariant toric metric
$g_0=\frac{1}{2\sqrt{-1}}\sum|\xi_\alpha|^{-2}d\xi_\alpha\wedge
d\bar{\xi}_\alpha$ on ${(\C^*)}^{n+1}$, in the neighborhood with
image sufficiently near $R$ under the map $\Log_t$.

Next we works on a neighborhood of an vertex $x\in C^{\wedge}_{i}$.
Near this point, there are $(n+1)$ terms
$$
t^{-v(j^{(1)})}z^{j^{(1)}},\, \ldots, t^{-v(j^{(n+1)})}z^{j^{(n+1)}},
$$
of $f$ that are comparable to $t^{-v(i)}z^{i}$ and dominating other
terms as $t\to+\infty$. Denote $j^{(0)}=i$ and
$\zeta_p=t^{-v(j^{(p)})}z^{j^{(p)}}$ for $p=0,\ldots,n+1$. Then
$$
0=f(z)=\sum_{p=0}^{n+1} \zeta_p + \cdots,
$$
and
\begin{eqnarray*}
z_{\alpha}f_{z_\alpha}&=&\sum_j j_\alpha t^{-v(j)}z^j\\
&=&\sum_{p=0}^{n+1}j^{(p)}_{\alpha} \zeta_p + \cdots \\
&=& j^{(0)}_{\alpha}\left( \sum_{p=0}^{n+1} \zeta_p \right) +
\sum_{p=1}^{n+1}(j^{(p)}_{\alpha}-j^{(0)}_\alpha) \zeta_{p}
+ \cdots\\
&=&j^{(0)}_{\alpha}\left( f+\cdots \right) + \sum_{p=1}^{n+1}(j^{(p)}_{\alpha}-j^{(0)}_\alpha) \zeta_{p}
+ \cdots
\end{eqnarray*}
where ``$\cdots$" denotes the terms in $f$, up to multiple of a
constant, other than $t^{-v(j^{(p)})}z^{j^{(p)}}$, $p=0,\ldots,n+1$.
Note that for $\alpha=1,\ldots,n+1$,
$$
z_\alpha=t^{x_\alpha}b_{\alpha} +\cdots
$$
with certain leading coefficients $b=(b_1,\ldots,b_{n+1})$,  where
$x_\alpha$'s are the coordinates of the vertex $x\in C^{\wedge}_i$.
Hence
$$
\zeta_p=t^{l_{v, j^{(p)}}(x)}b^{j^{(p)}} +\cdots,
$$
where $l_{v, j^{(p)}}(x)=\langle x, j^{(p)}\rangle -v(j^{(p)})$.
We claim that for any $b$, there exists $\alpha\in\{1,\ldots,n+1\}$ such that
$$
\lim_{t\to+\infty}\sum_{p=1}^{n+1}(j^{(p)}_{\alpha}-j^{(0)}_\alpha)
\frac{\zeta_{p}}{\zeta_{0}}\neq 0.
$$
In fact, if it is not true, then by taking $t\to+\infty$, we have for all $\alpha$,
$$
\lim_{t\to+\infty}\sum_{p=1}^{n+1}(j^{(p)}_{\alpha}-j^{(0)}_\alpha)
\frac{\zeta_{p}}{\zeta_{0}}=\sum_{p=1}^{n+1}(j^{(p)}_{\alpha}-j^{(0)}_\alpha) \frac{b^{j^{(p)}}}{b^{j^{(0)}}} = 0.
$$
By using $0=f(z)=\sum\zeta_p +\cdots$, we also have
$$
1+ \sum_{p=1}^{n+1} \frac{b^{j^{(p)}}}{b^{j^{(0)}}} =0.
$$
Then, by comparing with the
polynomial $1+z_1+\cdots+z_{n+1}$, it is easy to see that
$M=(j^{(p)}_{\alpha}-j^{(0)}_\alpha)_{\alpha,p=1,\ldots,n+1}$
is the matrix of the affine transformation that maps the
neighborhood of the vertex $x$ to the primitive complex $\Sigma_n$.
Since $\Pi_v$ is a maximal dual complex, $M$ is invertible. This
implies $\frac{b^{j^{(p)}}}{b^{j^{(0)}}} =0$ for all
$p$ which is a contradiction.

Using the claim, we see that for large $t$ and each $z$ in the neighborhood, $f_{z_\alpha}\neq 0$ for some $\alpha$.
And
$$
\Omega_i\to
\frac{(-1)^{\alpha-1}b^{j^{(0)}}}{\sum_{p=1}^{n+1}(j^{(p)}_{\alpha}-j^{(0)}_\alpha)
b^{j^{(p)}}}
\frac{d|\xi_1|}{|\xi_1|}\wedge\cdots\widehat{\left(\frac{d|\xi_\alpha|}{|\xi_\alpha|}\right)}
\cdots\wedge\frac{d|\xi_{n+1}|}{|\xi_{n+1}|}.
$$
This shows that $\Omega_i$ is non-vanishing near a vertex  for large
$t$. Similarly, we can prove that $\Omega_i$ is non-vanishing near
any face with dimension between $1$ and $n$. This completes the
proof of the second statement.

Finally for the last statement of the theorem,  we observe that on
any compact subset $B\subset\CP^{n+1}\setminus {U}^{\wedge}_{i,t}$,
$t^{-v(i)}z^{i}$ is no longer a dominating term near $V_t\cap B$ and
hence $\Omega_i \to 0$ locally in $B$ as $t\to+\infty$ by the above
local expression of $\Omega_i$. This completes the proof of the
theorem.

\end{proof}

\noindent{\em Remark:} In his preprint \cite{MiNote},  Mikhalkin
defined a version of ``regular 1-form" on tropical curves. From the
proof of the theorem \ref{thm-form}, we see that in this case
($n=1$), for each $i\in\triangle_{d,\Z}^0$, the holomorphic $1$-form
$\Omega_{i,t}$ tends to a limit of the form

$$
 \frac{(-1)^{\alpha-1}}{i_{\alpha}-j^{(1)}_\alpha}
dx_1\wedge\cdots\wedge\widehat{dx_\alpha} \wedge\cdots\wedge dx_{2},
$$
where $x_\alpha=\displaystyle\lim_{t\to +\infty}
\Log_t|z|=\displaystyle\lim_{t\to +\infty} \Log|\xi|$. And this
limit is in fact a ``regular $1$-form" on the tropical variety
$\calA_K(V)=\Pi_v$ in the sense of Mikhalkin and dual to the
$1$-cycle
$C^{\wedge}_i=\displaystyle\lim_{t\to\infty}\Log_t(U^{\wedge}_{i,t})$.

Note that the  set of $n$-cycles $\{
C^{\wedge}_i\}_{i\in\triangle_{d,\Z}^0}$, constructed in the theorem
\ref{thm-form} does not cover $\Pi_v$. Hence
$\{{U}^{\wedge}_{i,t}\}_{i\in\triangle_{d,\Z}^0}$ cannot cover
$V_t$. In order to obtain an open covering
$\{U_{i,t}\}_{i\in\triangle_{d,\Z}^0}$ of $V_t$, we need to enlarge
${U}^{\wedge}_{i,t}$ suitably. In the case of our interest, we have
the following

\begin{thm}\label{thm-set}
Let $v:\triangle_{d,\Z}\to \R$ be a function satisfying the
condition of lemma \ref{lem-v},
$f_t=\displaystyle\sum_{j\in\triangle_{d,\Z}}t^{-v(j)}z^j$ ($t>0$)
be the patchworking polynomial of degree $d$ defined by $v$ with
non-Archimedean amoeba $\calA_K=\Pi_v$. Denote
$V_t=\{f_t=0\}\subset\CP^{n+1}$. Then for all $t>0$, there exists a
basis $\{\Omega_{i,t}\}_{i\in\triangle_{d,\Z}^0}$ of $H^{n,0}(V_t)$,
and open subsets ${U}_{i,t} \subset V_t$ such that
for each $i\in\triangle_{d,\Z}^0$,

\begin{enumerate}
\item $\Log_t({U}_{i,t})$ tends to an $n$-cycle $C_i$ such that $\Pi_v =\bigcup_{i\in\triangle_{d,\Z}^0}C_i$ and
$\{[C_i]\}_{i\in\triangle_{d,\Z}^0}$ forms a basis of
$H_n(\calA_K(V))$,

\item for any compact subset $B\subset \CP^{n+1}\setminus {U}_{i,t}$, $\Omega_{i,t}$
tends to zero in $V_t^o\cap B$ uniformly with respect to
$H^*_t(g_0)$.
\end{enumerate}
\end{thm}
\begin{proof}
By the pair-of-pants decomposition theorem \ref{thm-pop}, $V_t$ is a
union of pairs-of-pants $\overline{\calP}_n$'s. And the set of
pair-of-pants $\overline{\calP}_n$ are in one-one correspondence to
the set of vertices of $\Pi_v$. Hence each $\overline{\calP}_n$
corresponds to a simplex in the lattice subdivision of $\triangle_d$
and vice versa. Moreover, if the pair-of-pants
$\overline{\calP}_n(s)$ corresponds to the simplex $\sigma_s$ of the
lattice subdivision of $\triangle_d$, then $\overline{\calP}_n(s)$
is the closure of the preimage $\lambda^{-1}(\calU_s)$ of the
primitive piece $\calU_s$ dual to the simplex $\sigma_s$.

Now, for each $i\in\triangle_{d,\Z}^0$, conditions of lemma
\ref{lem-v} says that the lattice subdivision of $\triangle_d$
corresponding to $v$ restricts to a lattice subdivision of
$i-{\i}+\triangle_{n+2}$. Let $\{\sigma_s\}$ be the set of simplices
of the lattice subdivision of $i-{\i}+\triangle_{n+2}$. Let
$\Lambda_i$ be the set of $s$ such that for all
$j\in\triangle_{d,\Z}^0\setminus\{ i\}$, $j$ is not a vertex of
$\sigma_s$. We define
$$
C_i=C_i^{\wedge}\bigcup\left( \bigcup_{s\in \Lambda_i} \overline{\calU_{s}} \right),
$$
where $\calU_s$ is the primitive piece dual to the simplex $\sigma_s$.

\pagebreak
\begin{figure}[h]\label{Figure 7}
\centerline{\includegraphics[width=10cm]{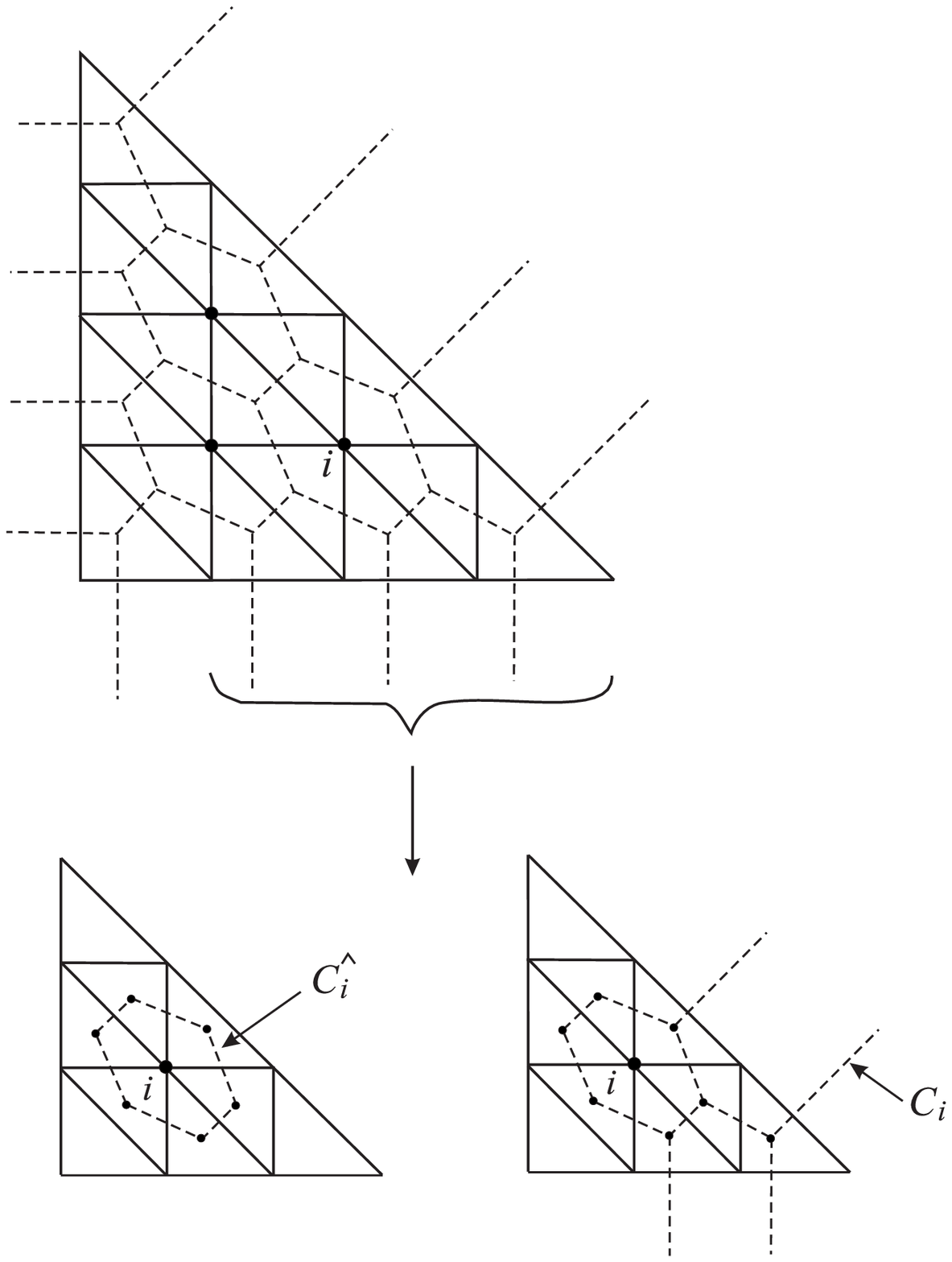}}
\end{figure}
\begin{center}
{\footnotesize Figure 7: Illustration of the sets $C^{\wedge}_{i}$
and $C_i$.}
\end{center}

Then we clearly have
$$
\Pi_v=\bigcup_{i\in\triangle_{d,\Z}^0} C_i
$$
and $[C_i]=[C^{\wedge}_i]$ in $H_n(\Pi_v)$. For these $C_i$, we define their $1/t$-neighborhood $C_{i,t}$ in $\Pi_v$ similar to those for $C^{\wedge}_{i,t}$ and define the preimages to be our enlarged open sets
$$
U_{i,t}=\lambda^{-1}\left( C_{i,t} \right).
$$
Then it is clear that $U_{i,t}\supset {U}^{\wedge}_{i,t}$ for all $i\in\triangle_{d,\Z}^0$ and
$$
V_t=\bigcup_{i\in\triangle_{d,\Z}^0} U_{i,t}.
$$

Finally, as $U^{\wedge}_{i,t}\subset U_{i,t}$, the second statement
of the theorem follows trivially from the corresponding statement in
theorem \ref{thm-form}. The proof is completed.
\end{proof}

\begin{figure}[h]\label{Figure 8}
\centerline{\includegraphics[width=8cm]{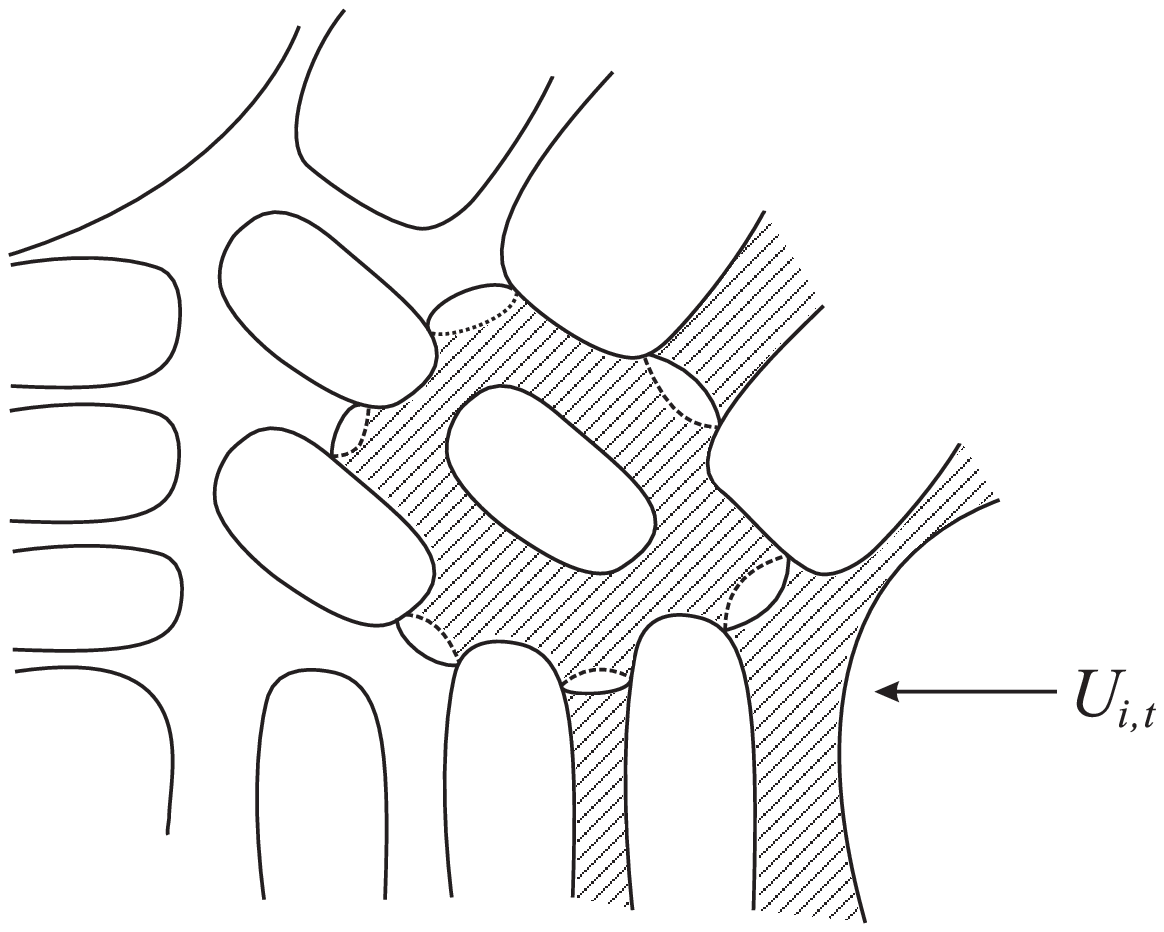}}
\end{figure}
\begin{center}
{\footnotesize Figure 8: Illustration of an open set $U_{i,t}$.}
\end{center}

\noindent{\bf Proof of the main theorem} By theorems \ref{thm-form}
and \ref{thm-set}, we remain to show the following
\begin{enumerate}
\item $U_{i,t}$ is close to an open subset of a
Calabi-Yau hypersurface $Y_{i,t}$ after normalization $H_t$,

\item $\Omega_{i,t}$ is close to a scalar multiple of the unique
holomorphic volume form $\Omega_{Y_{i,t}}$ of $Y_{i,t}$ on
$U_{i,t}$,

\item $\Omega_{i,t}$ is non-vanishing on the whole $U_{i,t}$.
\end{enumerate}
Note that the last item is needed because we have no information of
$\Omega_{i,t}$ at the points in $U_{i,t}\cap \left(
\CP^{n+1}\setminus(\C^*)^{n+1}\right)$ and the theorem
\ref{thm-form} shows that $\Omega_{i,t}$ tends to $0$ as
$t\to+\infty$ in $U_{i,t}\setminus U^{\wedge}_{i,t}$. We need to
show that even they are tending to $0$, $\Omega_{i,t}$ is still
non-vanishing on the whole $U_{i,t}$ for large but finite $t$.

To complete the proof of the main theorem, we take $v$ to be the
function given in the lemma \ref{lem-v},
$f_t=\displaystyle\sum_{j\in\triangle_{d,\Z}}t^{-v(j)}z^j$ be the
patchworking polynomial defined by $v$ and $V_t=\{ f_t=0\}$. Recall
that the lemma \ref{lem-v} implies that the subdivision
corresponding to $v$ restricts to a lattice subdivision of
$i-{\i}+\triangle_{n+2}$. Consider the truncated polynomial
$f_{i,t}=\displaystyle\sum_{j \in i-{\i}+\triangle_{n+2}}
t^{-v(j)}z^j$. This truncated polynomial $f_{i,t}$ factorizes as
$$
f_{i,t}= z^{i-{\i}} \sum_{j \in \triangle_{n+2}} t^{-v(j+i-{\i})}z^j
$$
So the variety $\{f_{i,t}=0\}\cap(\C^*)^{n+1}$ can be regarded as an
open subset of the Calabi-Yau variety $Y_{i,t}$ defined by the polynomial
${f}^{CY}_{i,t}=\displaystyle\sum_{j \in \triangle_{n+2}} t^{-v(j+i-{\i})}z^j$.

In \cite{Mi03}, it was shown that the normalized varieties $H_t(V_t)$
converge in the Hausdroff metric to the lift $W(\calA_K)$ of the
corresponding non-Archimedean amoeba $\calA_K$, where the lift
$W(\calA_K)$ is the image of $\calA_K$ under the map
$W:(K^*)^{n+1}\to (\C^*)^{n+1}$ defined as
$$
W(b_1,\ldots,b_{n+1})=\left( e^{-val_K(b_1)+i\arg(b_1^m)},\ldots,
e^{-val_K(b_{n+1})+i\arg(b_{n+1}^m)} \right),
$$
where $b_\alpha^m\in\C$ is the coefficient of
$t^{-val_K(b_\alpha)}$, i.e, the leading coefficient of $b_\alpha$.
For simplicity, we call $\arg(b_1^m),\ldots,\arg(b_{n+1}^m)$ the leading arguments of $b\in
K$. It was shown that $W(\calA_K)$ depends only on the leading
arguments of the coefficients of the defining polynomial $f$.

We can define similarly the lift $W(C_{i,t})$ for each $i\in\triangle_{d,\Z}^0$. Then the proof in
\cite{Mi03} applies directly to show that $W(C_{i,t})$ depends only
on the leading arguments of the coefficients of the terms that
determines $C_{i,t}$ and the normalized open set $H_t(U_{i,t})$ is
close to $W(C_{i,t})$. Since $f_t$ and the truncated polynomial
$f_{i,t}$ have the same coefficients corresponding to $C_{i,t}$, we
see immediately that $H_t(U_{i,t})$ is close to $H_t({U}^{CY}_{i,t})$
in Hausdroff distance, where ${U}^{CY}_{i,t}$ is the preimage of the
neighborhood ${C}^{CY}_{i,t}$ corresponding to $C_{i,t}$ in the lattice
subdivision of $i-{\i}+\triangle_{n+2}$. Therefore, $U_{i,t}$ is
close to an open set of a Calabi-Yau manifold $Y_{i,t}$ after normalization
$H_t$. This proves the first statement.

To see the other two statements, we observe that the limiting
behavior of $\Omega_{i,t}$ shows that $\Omega_{i,t}$ is close to the
corresponding holomorphic $n$-form $\Omega_{Y_{i,t}}$ of the
Calabi-Yau hypersurface
$\{{f}^{CY}_{i,t}=0\}=Y_{i,t}\subset\CP^{n+1}$. In fact,  since
${f}^{CY}_{i,t}=0$ on the hypersurface $\{f_{i,t}=0\}$, we have
\begin{eqnarray*}
\Omega_{Y_{i,t}} & = & (\log t)^{-n} \left[
\frac{(-1)^{\alpha-1}t^{-v(i)}z^{\i}}{z_\alpha(
{f}^{CY}_{i,t})_{z_{\alpha}}} \right]
\frac{dz_1}{z_1}\wedge\cdots\widehat{\left(\frac{dz_\alpha}{z_\alpha}\right)}
\cdots\wedge\frac{dz_{n+1}}{z_{n+1}} \\
 & = & (\log t)^{-n}  \left[ \frac{(-1)^{\alpha-1}t^{-v(i)}z^{i}}{z_\alpha z^{i-{\i}}( {f}^{CY}_{i,t})_{z_{\alpha}}} \right]
\frac{dz_1}{z_1}\wedge\cdots\widehat{\left(\frac{dz_\alpha}{z_\alpha}\right)}
\cdots\wedge\frac{dz_{n+1}}{z_{n+1}}\\
 & = & (\log t)^{-n}
\left[ \frac{(-1)^{\alpha-1}t^{-v(i)}z^{i}}{z_\alpha (z^{i-{\i}} {f}^{CY}_{i,t})_{z_{\alpha}}} \right]
\frac{dz_1}{z_1}\wedge\cdots\widehat{\left(\frac{dz_\alpha}{z_\alpha}\right)}
\cdots\wedge\frac{dz_{n+1}}{z_{n+1}} \\
& = &
(\log t)^{-n} \left[
\frac{(-1)^{\alpha-1}t^{-v(i)}z^{i}}{z_\alpha (f_{i,t})_{z_{\alpha}}} \right]
\frac{dz_1}{z_1}\wedge\cdots\widehat{\left(\frac{dz_\alpha}{z_\alpha}\right)}
\cdots\wedge\frac{dz_{n+1}}{z_{n+1}}.
\end{eqnarray*}
Therefore using the fact that $f$ and $f_{i,t}$ contain the same
dominating terms on $U_{i,t}$, we see that $\Omega_{i,t}$ is close
to $\Omega_{Y_{i,t}}$ in the sense described in theorem
\ref{thm-form}.

In this local coordinate, the ratio
$$
\frac{\Omega_{i,t}}{\Omega_{Y_{i,t}}}=\frac{z_\alpha (f_{i,t})_{z_{\alpha}}}{z_\alpha f_{z_{\alpha}}}
$$
is a holomorphic function. Note that $f_{i,t}$ and $f$ contain the
same dominating terms in the neighborhoods corresponding to faces of
$C_{i,t}$ and ${C}^{CY}_{i,t}$ respectively. In a neighborhood of
$U_{i,t}\cap\left(\CP^{n+1}\setminus(\C^*)^{n+1}\right)$, we
consider those open subsets correspond to top dimensional faces. In
these open subsets, the polynomial $f$ and $f_{i,t}$ are dominated
by exactly two terms $t^{j^{1}}z^{j^{(1)}}$ and
$t^{j^{2}}z^{j^{(2)}}$. Then
\begin{eqnarray*}
\frac{\Omega_{i,t}}{\Omega_{Y_{i,t}}} &=& \frac{z_\alpha (f_{i,t})_{z_{\alpha}}}{z_\alpha f_{z_{\alpha}}} \\
&=& \frac{z_\alpha \left( t^{j^{1}}z^{j^{(1)}}+ t^{j^{2}}z^{j^{(2)}}
+\cdots \right)_{z_\alpha}}{z_\alpha \left( t^{j^{1}}z^{j^{(1)}}+ t^{j^{2}}z^{j^{(2)}}
+\cdots \right)_{z_\alpha} +z_\alpha \left( f-f_{i,t} \right)_{z_\alpha} } \\
&=& \frac{j^{(1)}_{\alpha}t^{j^{1}}z^{j^{(1)}}+ j^{(2)}_{\alpha}t^{j^{2}}z^{j^{(2)}}
+\cdots  }{ j^{(1)}_{\alpha}t^{j^{1}}z^{j^{(1)}}+ j^{(2)}_{\alpha}t^{j^{2}}z^{j^{(2)}} +\cdots } \\
&=& \frac{\left( j^{(1)}_{\alpha}- j^{(2)}_{\alpha}\right)
t^{j^{1}}z^{j^{(1)}}+ j^{(2)}_{\alpha} \left(f_{i,t}+\cdots\right)
+\cdots  }{ \left( j^{(1)}_{\alpha}- j^{(2)}_{\alpha}\right)
t^{j^{1}}z^{j^{(1)}}+ j^{(2)}_{\alpha} \left(f +\cdots\right)
+\cdots }.
\end{eqnarray*}
As before, ``$\cdots$" means a linear combination of the terms other
than $t^{j^{1}}z^{j^{(1)}}$ and $t^{j^{2}}z^{j^{(2)}}$. Therefore
${\Omega_{i,t}}/{\Omega_{Y_{i,t}}}$ is non-vanishing in these open
subsets. Taking closure of these open subsets in the neighborhood,
we conclude that ${\Omega_{i,t}}/{\Omega_{Y_{i,t}}}$ is
non-vanishing on
$U_{i,t}\cap\left(\CP^{n+1}\setminus(\C^*)^{n+1}\right)$. Since
$\Omega_{Y_{i,t}}$ is the holomorphic volume of the Calabi-Yau
hypersurface $Y_{i,t}$, it is non-vanishing and hence $\Omega_{i,t}$
is also non-vanishing on the whole $U_{i,t}$.

Finally, as $U_{i,t}\supset {U}^{\wedge}_{i,t}$, the last statement follows immediately from the last statement of the theorem \ref{thm-form}
This completes the proof of the main theorem.

\subsection{Asymptotically special Lagrangian fibers}\label{sec-asl}
From the fibration $\lambda$ given in \cite{Mi03}, for each $n$-cell
$e$ of $\Pi_v$, there exists a point $x\in e$
such that the fiber $\lambda^{-1}(x)$ is a Lagrangian $n$-torus
$\T^n\subset V$ which is actually given by $\{z\in V_t: \Log_t|z|=x\}$.
Therefore, when restricted to this fiber
$$
\left.\Omega_i\right|_{\lambda^{-1}(x)} = \left(\frac{\sqrt{-1}}{\log t}\right)^{n} \left[
\frac{(-1)^{\alpha-1}t^{-v(i)}z^{i}}{(i_{\alpha}-j^{(1)}_\alpha)
t^{-v(i)}z^{i}+ \cdots} \right]
d\theta_1\wedge\cdots\widehat{d\theta_\alpha}
\cdots\wedge d\theta_{n+1}.
$$
where $\theta_\alpha=\arg z_\alpha$ for $\alpha=1,\ldots, n+1$.
So we have
\begin{thm}\label{thm-fibration}
Let $V_t$ be the family of smooth hypersurfaces in $\CP^{n+1}$ of
degree $d$, $U_{i,t}$ be the open sets and $\Omega_{i,t}$ be the
holomorphic $n$-forms in the main theorem. Let $\lambda: V_t \to
\Pi_v$ be the stratified $\T^n$-fibration given by Mikhalkin
\cite{Mi03}. Then for any $i\in\triangle_{d,\Z}^0$ and any $n$-cell
$e$ of $C^{\wedge}_{i,t}$ used in the proof of the main theorem,
there exists $x\in e$, independent of $t$, such that for all
$j\in\triangle_{d,\Z}^0 $,
$$
\lim_{t\to+\infty}\operatorname{Im}\left.\left(e^{\frac{n\pi\sqrt{-1}}{2}}\Omega_{j,t}
\right)\right|_{\lambda^{-1}(x)} = 0 \quad\mbox{and}
$$
$$\lim_{t\to+\infty}\operatorname{Re}\left.\left(e^{\frac{n\pi\sqrt{-1}}{2}}\Omega_{j,t}
 \right)\right|_{\lambda^{-1}(x)} \quad\mbox{is non-vanishing}.
$$
\end{thm}
In particular, the Lagrangian fibers given by Mikhalkin are in fact
``asymptotically special Lagrangian of phase $n\pi/2$" with respect to the holomorphic $n$-form $\Omega_{i,t}$
constructed in the main theorem.

\subsection{Hypersurfaces in other toric varieties; the case of curves}\label{sec-2d}
It is clear from the works of Mikhalkin \cite{Mi03}, our result can
be modified to include other toric varieties such as
$\CP^m\times\CP^n$. In particular, if we apply our method to curves
in $\CP^1\times\CP^1$ instead of $\CP^2$, we will obtain stronger
results for Riemann surfaces. It applies to Riemann surfaces of all
genus, not just for $g= \frac{(d-1)(d-2)}{2}$. And the Calabi-Yau
components  actually form a connected sum decomposition of the
curve. This is a fact which is probably not true in higher
dimensions. The key issue is that in this case, one can obtain a
subdivision with dual complex $\Pi_v$ with 1-cycles $\{C_i\}$ such
that $C_i\cap C_j=\emptyset$ for $i\neq j$. In summary, we have
\begin{thm}\label{thm-main-2d}
For any integer $g\ge 1$, there is a family of  smooth genus g
curves $V_t$ of bi-degree $(g+1,2)$ in $\CP^1\times\CP^1$, such that
$V_t$ can be written as
$$
V_t=\bigcup_{i=1}^g U_{i,t}
$$
where $\{ U_{i,t}\}$ is a family of closed subsets $U_{i,t}\subset
V_t$ such that topologically $V_t = U_{1,t}\# \cdots \# U_{g,t}$,
the connected sum of $U_{i,t}$, $i=1,\ldots, g$ and after
normalization $H_t:(\C^*)^2 \to (\C^*)^2$ defined by
$$
H_t(z_1,z_2)=\left( |z_1|^{\frac{1}{\log t}}\frac{z_1}{|z_1|} ,
|z_2|^{\frac{1}{\log t}}\frac{z_2}{|z_2|} \right),
$$
\begin{enumerate}
\item
$U_{i,t}$ is close in Hausdroff distance on
$(\C^*)^2$ to an open subset of an elliptic
curves $Y_{i,t}$ in $\CP^1\times\CP^1$;
\item there exists a basis $\{\Omega_{i,t}\}_{i=1}^g$ of $H^{1,0}(V_t)$
such that for each $i=1,\ldots ,g$, $\Omega_{i,t}$ is nonvanishing
and $\epsilon$-closed to the holomorphic $1$-form $\Omega_{Y_{i,t}}$
of $Y_{i,t}$ on $U_{i,t}$ with respect to the metric induced from
the pull-back $H^*_t(g_0)$;

\item for any compact subset
$B\subset (\C^*)^2\setminus U_{i,t}\subset  \CP^1\times\CP^1$,
$\Omega_{i,t}$ tends to zero in $V_t\cap B$ uniformly with respect
to $H^*_t(g_0)$.
\end{enumerate}
\end{thm}
\begin{proof}
All the steps in the proof of the main theorem apply and give all
results concerning the family $U_{i,t}$ except the assertion about
the connected-sum. To see this, we need the existence of certain
maximal lattice subdivision of the Newton polytope
$\triangle=[0,g+1]\times [0,2]$ of a generic curves of bi-degree
$(g+1,2)$ such that the set of vertices $\mbox{ver}(\sigma)$ of each
simplex $\sigma$ contains at most one interior lattice points of
$\triangle$. For this kind of subdivisions, each primitive piece
associated to at most one interior lattice point. Then any two
cycles $C^{\wedge}_i$ and $C^{\wedge}_j$ must have empty
intersection for $i\neq j$. As demonstrated in the following figure,
the existence of such subdivision is easy to see. In fact, one shows
that there exists a function $v:\triangle_{\Z}\to\R$ such that
$\Pi_v$ is a maximal dual complex of $\triangle$ giving the required
subdivision.

\pagebreak
\begin{figure}[h]\label{Figure 9}
\centerline{\includegraphics[width=10cm]{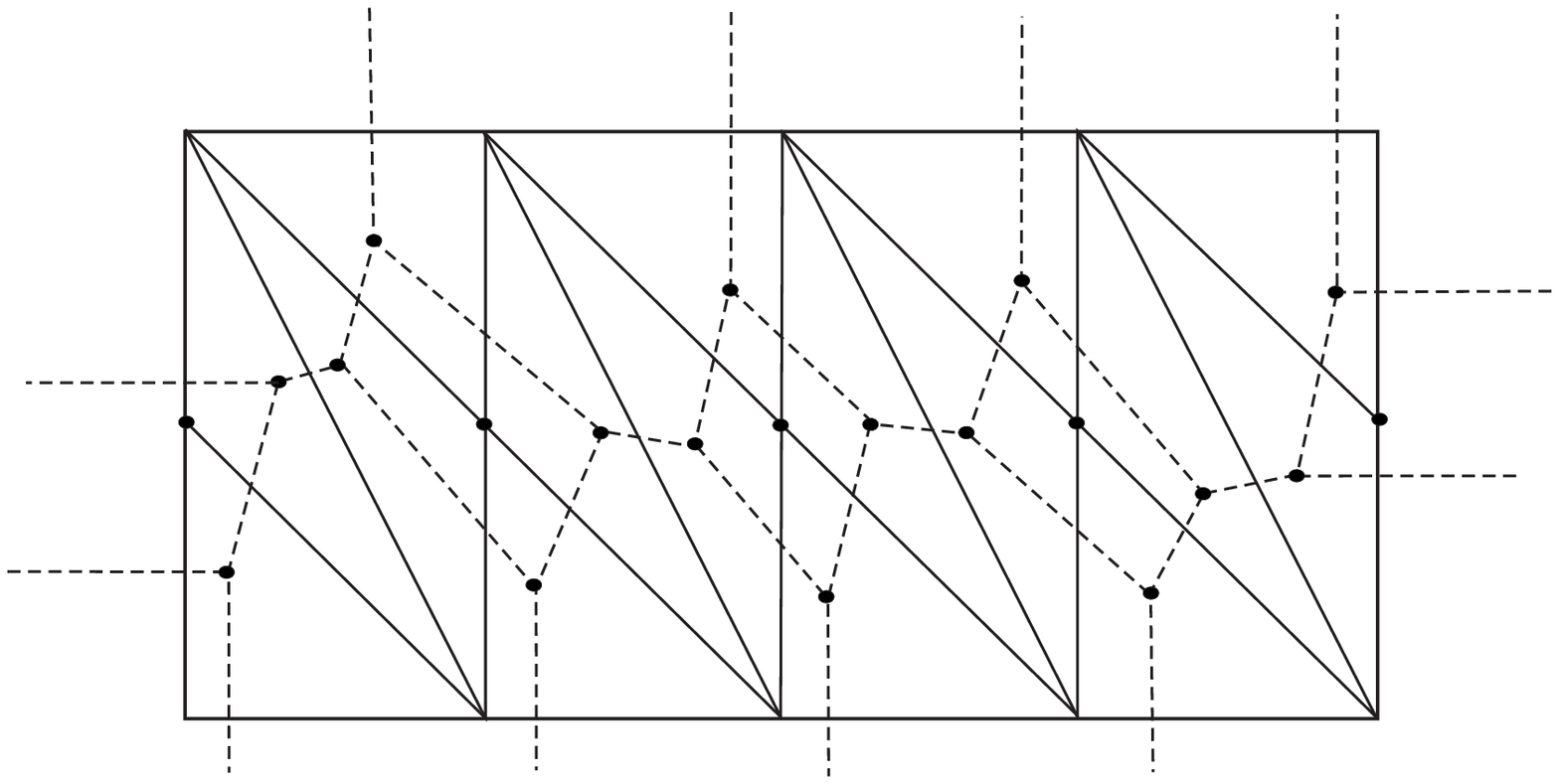}}
\end{figure}
\begin{center}
{\footnotesize Figure 9: A maximal dual $\triangle$ complex with the required properties.}
\end{center}

As in the proof of the main theorem, we can construct open  subsets
$U_{i,t}$ and holomorphic 1-form $\Omega_{i,t}$ corresponding to
each cycles $C_i$, where $i=1,\ldots,g$ is in one-one corresponding
with the set of interior lattice points of $\triangle=[0,g+1]\times
[0,2]$. This gives all the assertions except the part of the
connected sum.

Note that in this case, two cycles $C^{\wedge}_i$ and $C^{\wedge}_j$
of $\Pi_v$ corresponding to different interior lattice points of
$\triangle$ do not intersect. In fact, they are, at least, separated
by an  edge $e$ of $\Pi_v$. Therefore, either $C_i\cap
C_j=\emptyset$ or $C_i\cap C_j=\{x\}$ where $x\in e$ is the common
boundary of the corresponding primitive pieces. It is now clear that
$V_t$ are just all $U_{i,t}$ gluing along the circles
$\lambda^{-1}(x)$ of these $x$'s. So $V_t = U_{1,t}\# \cdots \#
U_{g,t}$ and proof is completed.
\end{proof}

We would like to remark that in this complex one dimensional case,
special Lagrangian submanifolds of $\Omega_{i,t}$ always exists.
Namely, they are given by the horizontal or vertical trajectories of
the quadratic differential $\Omega_{i,t}^2$. Therefore, we have
\begin{thm}\label{thm-fibration-2d}
Let $V_t$ be the family of smooth curves in $\CP^1\times\CP^1$ of
degree $g+1$, $\Omega_{i,t}$ be the
holomorphic $1$-forms in the theorem \ref{thm-main-2d}. Then for any
$i=1,\ldots,g$, there exists special Lagrangain fibration of
$\Omega_{i,t}$.
\end{thm}


\section{Proof of the key lemma}\label{sec-lem}
In this section, we prove the lemma \ref{lem-v} concerning the dual
complex $\Pi_v$ given by the function $v:\triangle_{d,\Z}\to
\R$ defined by
$$
v(j)=\sum_{\alpha=1}^{n+1} j_\alpha^2+\left( \sum_{\alpha=1}^{n+1}
j_\alpha \right)^2
$$
for $j=(j_1,\ldots,j_{n+1})\in \triangle_{d,\Z}$. To
simplify notation, we will write $j_0=\displaystyle\sum_{\alpha=1}^{n+1}
j_\alpha$. Then $v(j)=\displaystyle\sum_{\alpha=0}^{n+1} j_\alpha^2$. We start
with
\begin{lem}\label{lem-v1}
Let $i\in  \triangle_{d,\Z}$ and $j\in
\Z^{n+1}\setminus\left(\{0\}\cup\{\pm
e_\alpha\}_{\alpha=1}^{n+1}\cup \{e_\alpha-e_\beta
\}_{\alpha\neq\beta}\right)$ such that $i+j\in
\triangle_{d,\Z}$. Then there exists no $x\in\R^{n+1}$ such
that
\begin{eqnarray*}
\langle i, x \rangle -v(i) &=& \langle i+j, x \rangle -v(i+j) \\
 & > & \langle r, x \rangle -v(r),\quad r\in \triangle_{d,\Z}\setminus\{i,i+j\}.
\end{eqnarray*}
\end{lem}
\begin{proof}
We separate the proof into several steps.

\noindent{\bf Step 1:} If $i$ is an interior lattice point of
$\triangle_{d}$, then for all $x\in\R^{n+1}$ satisfying the
condition in the lemma, we have
$$
\begin{cases}
|x_\alpha-2(i_\alpha+i_0)|<2,& \forall\,
\alpha=1,\ldots, n+1\\
|x_\alpha-x_\beta-2(i_\alpha-i_\beta)|<2,& \forall\, \alpha\neq
\beta=1,\ldots, n+1.
\end{cases}
$$

\noindent{\em Proof of Step 1:} Since $i\in \triangle_{d,\Z}^0$, the
points $i\pm e_\alpha$, $i+e_\beta-e_\alpha$ belong to
$\triangle_{d,\Z}$, for all $\alpha$ and $\beta\neq
\alpha$. By assumption, $j\neq \pm e_\alpha$, $e_i-e_j$, we have the
following strict inequalities from the condition of the lemma
$$
\begin{cases}
\langle i\pm e_\alpha,x \rangle -v(i\pm e_\alpha) &< \langle i,x
\rangle -v(i)
\\
\langle i+e_\beta-e_\alpha,x \rangle -v(i+e_\beta-e_\alpha) &<
\langle i,x \rangle -v(i).
\end{cases}
$$
Using the definition of $v$, we have
$$
\begin{cases}
\pm x_\alpha - \left[ (\pm 2i_\alpha +1)+(\pm 2i_0 +1) \right] &< 0
\\
x_\beta-x_\alpha - \left[ ( 2i_\beta +1)+(-2i_\alpha +1) \right]  &<
0.
\end{cases}
$$
Interchanging the $\beta$ and $\alpha$ gives the required
inequalities
$$
\begin{cases}
|x_\alpha-2(i_\alpha+i_0)|<2,& \forall\, \alpha=1,\ldots, n+1\\
|x_\beta-x_\alpha-2(i_\beta-i_\alpha)|<2,& \forall\, \beta\neq
\alpha=1,\ldots, n+1.
\end{cases}
$$

\noindent{\bf Step 2:} Either $i$ or $i+j$ belong to boundary of
$\triangle_{d}$.

\noindent{\em Proof of Step 2:} Suppose not, then both $i$ and $i+j
\in \triangle_{d,\Z}^0$. Applying Step 1 to $i$
and $i+j$, we have the following inequalities
$$
\begin{cases}
|x_\alpha-2(i_\alpha+i_0)|<2,& \forall\, \alpha=1,\ldots, n+1\\
|x_\beta-x_\alpha-2(i_\beta-i_\alpha)|<2,& \forall\,\beta\neq
\alpha=1,\ldots, n+1.
\end{cases}
$$
and
$$
\begin{cases}
|x_\alpha-2(i_\alpha+j_\alpha+i_0+j_0)|<2,& \forall\, \alpha=1,\ldots, n+1\\
|x_\beta-x_\alpha-2(i_\beta+j_\beta-i_\alpha-j_\alpha)|<2,&
\forall\,\beta\neq \alpha=1,\ldots, n+1.
\end{cases}
$$
Therefore,
$$
\begin{cases}
|j_\alpha+j_0|<2,& \forall\,\alpha=1,\ldots, n+1\\
|j_\beta-j_\alpha|<2,& \forall\, \beta\neq\alpha=1,\ldots, n+1.
\end{cases}
$$
As $j_\beta$ are integers, the first inequality above implies for
all $\alpha=1,\ldots,n+1$,
$$
-1 \le j_\alpha+j_0 \le 1.
$$
Recalling $j_0=\displaystyle\sum_{\beta=1}^{n+1}j_\beta$ and summing over
$\alpha$, we have
$$
-(n+1)\le (n+2)j_0 \le n+1.
$$
Hence $j_0=\displaystyle\sum_{\beta=1}^{n+1}j_\beta=0$. On the other hand, the
second inequality above implies that either $j_\beta\ge 0$ for all
$\beta=1,\ldots,n+1$ or $j_\beta\le 0$ for all $\beta=1,\ldots,n+1$.
Therefore, we must have $j=0$ which is a contradiction. Step 2 is
proved.

\noindent{\bf Step 3:} Both $i$ and $i+j$ belong to the boundary of
$\triangle_d$.

\noindent{\em Proof of Step 3:} Suppose this is not true. We may
assume $i$ belongs to the interior and $i+j$ belongs to the
boundary. Then we can apply Step 1 to $i$ and get
$$
\begin{cases}
|x_\alpha-2(i_\alpha+i_0)|<2,& \forall\, \alpha=1,\ldots, n+1\\
|x_\beta-x_\alpha-2(i_\beta-i_\alpha)|<2,& \forall\, \beta\neq
\alpha=1,\ldots, n+1.
\end{cases}
$$
However, we do not have all the inequalities as at least one of the
points $i+j\pm e_\alpha$, $i+j+e_\beta-e_\alpha$ lies outside
$\triangle_d$.

\noindent{\em Case 1:} $i_0+j_0=\displaystyle\sum_{\beta=1}^{n+1}
(i_\beta+j_\beta) \le d-1$.

In this case, we still have
$$
\begin{cases}
|x_\alpha-2(i_\alpha+j_\alpha+i_0+j_0)|<2,& \forall\,
\alpha\mbox{ with } i_\alpha+j_\alpha\ge 1\\
|x_\beta-x_\alpha-2(i_\beta+j_\beta-i_\alpha-j_\alpha)|<2,&
\forall\,\beta\neq\alpha\mbox{ with } i_\beta+j_\beta,\,
i_\alpha+j_\alpha\ge 1.
\end{cases}
$$
And for those $\alpha$ with $i_\alpha+j_\alpha=0$,  we only have the
one-sided inequality
$$
x_\alpha-2(i_\alpha+j_\alpha+i_0+j_0)<2.
$$
Hence for $\alpha$ with $i_\alpha+j_\alpha\ge 1$, we still have
$$
|j_\alpha+j_0|<2.
$$
However, for $\alpha$ with $i_\alpha+j_\alpha=0$, we only have
$$
-2<i_\alpha+j_0.
$$
Since $i+j\in\partial\triangle_d$ and we are assuming
$i_0+j_0=\displaystyle\sum_{\beta=1}^{n+1} (i_\beta+j_\beta) \le d-1$ in this
case, there exists $\alpha_o$ such that
$i_{\alpha_o}+j_{\alpha_o}=0$. For this $\alpha_o$, we get
$$
-2<j_{\alpha_o}+j_0=-i_{\alpha_o}+j_0.
$$
As $i\in\triangle_{d,\Z}^0$, we have $i_{\alpha_o}\ge 1$, and hence
$j_0\ge 0$.

If $i+j\neq 0$, then there also exists $\alpha$ such that
$i_\alpha+j_\alpha\ge 1$. For all these $\alpha$,
$i+j+e_{\alpha_o}-e_\alpha \in\triangle_{d,\Z}$ and we can
apply the condition of the lemma to get the strict inequality
$$
\langle i+j+e_{\alpha_o}-e_\alpha,x \rangle
-v(i+j+e_{\alpha_o}-e_\alpha) < \langle i+j, x \rangle -v(i+j).
$$
Using $i_{\alpha_o}+j_{\alpha_o}=0$, this gives
$$
x_{\alpha_o}-x_\alpha+2(i_\alpha+j_\alpha)<2.
$$
Together with $|x_{\alpha_o}-x_\alpha-2(i_{\alpha_o}-i_\alpha)|<2$,
we arrive at
$$
-2<2\left( 1-(i_\alpha+j_\alpha) \right) -2(i_{\alpha_o}-i_\alpha),
$$
which implies
$$
j_\alpha<2-i_{\alpha_o}\le 1.
$$
Therefore, $i_\alpha\le 0$ for those $\alpha$ with
$i_\alpha+j_\alpha\ge 1$. Putting these into $0\le
j_0=\displaystyle\sum_{\beta=1}^{n+1}j_\beta$, we have
$$
0\le \sum_{\{\alpha_o\,:\,
i_{\alpha_o}+j_{\alpha_o}=0\}}j_{\alpha_o} = - \sum_{\{\alpha_o\,:\,
i_{\alpha_o}+j_{\alpha_o}=0\}}i_{\alpha_o} <0,
$$
which is a contradiction. So we must have $i+j=0$, that is
$i_\alpha+j_\alpha=0$ for all $\alpha$. Then for all $\alpha$,
$i+j+e_\alpha=e_\alpha\in\triangle_{d,\Z}$. The condition of
the lemma implies
$$
\langle i+j+e_\alpha, x \rangle -v(i+j+e_\alpha)<\langle i+j,
x\rangle -v(i+j).
$$
So $x_\alpha<v(e_\alpha)=2$ for all $\alpha$. On the other hand, the
condition of the lemma gives the equality
$$
\langle i+j, x \rangle -v(i+j)= \langle i, x\rangle -v(i),
$$
which is
$$
\sum_{\beta=0}^{n+1}i_\beta^2=\langle i, x\rangle.
$$
Using $x_\alpha<2$ and $i_\alpha\ge 1$, we have
$$
\sum_{\beta=0}^{n+1}i_\beta^2 < 2 \sum_{\beta=1}^{n+1}i_\beta=2i_0.
$$
Hence,
$$
\sum_{\beta=1}^{n+1}i_\beta^2 + \left( i_0 -1 \right)^2 < 1.
$$
This is a contradiction as $i_\alpha\ge 1$ for all $\alpha$. So we
have proved that the case 1 with assumption $i_0+j_0\le d-1$ is
impossible and hence we must be in the situation of the following

\noindent{\em Case 2:} $i_0+j_0=d$.

In this case, for any $\alpha$ with $i_\alpha+j_\alpha\ge 1$, we
have $i+j-e_\alpha\in\triangle_{d,\Z}$. Therefore, the
following strict inequality is satisfied
$$
\langle i+j-e_\alpha, x\rangle -v(i+j-e_\alpha) < \langle i+j,
x\rangle -v(i+j).
$$
Using $i_0+j_0=d$, we get
$$
-2< x_\alpha-2(i_\alpha+j_\alpha+d)
$$
provided $i_\alpha+j_\alpha\ge 1$. One also have
$i+j+e_\gamma-e_\alpha\in\triangle_{d,\Z}$ for any
$\gamma\neq \alpha$. Similar argument implies
$$
x_\gamma-x_\alpha-2(i_\gamma+j_\gamma-i_\alpha-j_\alpha)<2
$$
provided $i_\alpha+j_\alpha\ge 1$ and $\gamma\neq \alpha$.

The second inequality together with
$$
|x_\gamma-x_\alpha-2(i_\gamma-i_\alpha)|<2
$$
imply
$$
-2<x_\gamma-x_\alpha-2(i_\gamma-i_\alpha)<2+2(j_\gamma-j_\alpha).
$$
That is
$$
j_\alpha<2+j_\gamma
$$
for $\alpha$ with $i_\alpha+j_\alpha\ge 1$ and $\gamma\neq \alpha$.

Suppose there is an $\gamma_o$ such that
$i_{\gamma_o}+j_{\gamma_o}=0$, then $j_{\gamma_o}\neq j$ and hence
$$
j_\alpha<2+j_{\gamma_o}=2-i_{\gamma_o} <1
$$
since $i\in\triangle_{d,\Z}^0$. So $j_\alpha\le 0$
for any $\alpha$ with $i_\alpha+j_\alpha\ge 1$. It is trivial that
$j_\alpha\le 0$ for those $\alpha$ with $i_\alpha+j_\alpha=0$, we
have $j_\alpha\le 0$ for all $\alpha$ which in turn implies that
$i_0=d-j_0=d-\displaystyle\sum_{\beta=1}^{n+1}j_\beta \ge d$. This is impossible
as $i$ is in the interior of $\triangle_d$. Therefore, we cannot
have $\gamma_o$ such that $i_{\gamma_o}+j_{\gamma_o}=0$. That is,
$i_\alpha+j_\alpha\ge 1$ for all $\alpha$.

Now the other two inequalities must be satisfied, i.e.
$$
|x_\alpha-2(i_\alpha+i_0)|<2;
$$
$$
-2< x_\alpha-2(i_\alpha+j_\alpha+d).
$$
Combining these, we have
$$
-2<x_\alpha-2i_\alpha-2j_\alpha-2d<2+2i_0-2j_\alpha-2d.
$$
Hence
$$
j_\alpha+d< 2+ i_0 =2 +d-j_0,
$$
i.e.
$$
j_\alpha+j_0 <2.
$$
Summing over $\alpha$, we have $(n+2)j_0 <2(n+1)$ which implies
$j_0\le 1$. As $i_0+j_0=d$ and $i_0\le d-1$, we see that $j_0=1$
(and $i_0=d-1$). Putting it back into the inequality, we have
$j_\alpha<1$, hence, $j_\alpha\le 0$ for all $\alpha$. Summing over
to get the contradiction that $j_0\le 0$. This proved the Case 2 and
the proof of Step 3 is completed.

\noindent{\bf Step 4:} Either $j\in
\{y\in\R^{n+1}:y_\beta=0\}\cap\partial\triangle_d\cap\Z^{n+1}$ for
some $\beta$ and $i+j \in
\{y\in\R^{n+1}:y_\alpha=0\}\cap\partial\triangle_d\cap\Z^{n+1}$ for
some $\alpha$.

\noindent{\em Proof of Step 4:} Since Step 3 shows that both $i$ and
$i+j$ belong to $\partial\triangle_d\cap\Z^{n+1}$, if Step 4 is not
true, then either $i$ or $i+j$ belongs to the interior of the face
$\{y\in\R^{n+1}:\sum y_\beta=d\}\cap\partial\triangle_d$.

Let first assume that both $i$ and $i+j$ belong to the interior of
the face $\{y\in\R^{n+1}:\sum y_\beta=d\}\cap\partial\triangle_d$.
Then both $i+e_\beta-e_\alpha$ and $i+j+e_\beta-e_\alpha$ belong to
$\triangle_{d,\Z}$, we have the inequalities
$$
|x_\beta-x_\alpha-2(i_\beta-i_\alpha)|<2,
$$
$$
|x_\beta-x_\alpha-2(i_\beta+j_\beta-i_\alpha-j_\alpha)|<2.
$$
Hence $|j_\beta-j_\alpha|<2$. This implies $j_\alpha\ge 0$ for all
$\alpha$ or $j_\alpha\le 0$ for all $\alpha$. Together with
$$
j_0=\sum_{\beta=1}^{n+1}j_\beta=\sum_{\beta=1}^{n+1}(i_\beta+j_\beta)
-\sum_{\beta=1}^{n+1}i_\beta=d-d=0,
$$
we have $j=0$ which is a contradiction.

Now we may assume that $i+j$ belongs to the interior of the face
$\{y\in\R^{n+1}:\sum y_\beta=d\}\cap\partial\triangle_d$ but
$p\in\{y\in\R^{n+1}: y_\beta=0\}\cap\partial\triangle_d$ for some
$\beta$. Then we have
$$
\begin{cases}
i_\beta=0,& \\
i_\alpha+j_\alpha\ge 1,& \forall\,\alpha\\
i_0+j_0=0.&
\end{cases}
$$
If $i_0\le d-1$. Then $i+e_\alpha\in\triangle_{d,\Z}$ and
hence
$$
x_\alpha-2(i_\alpha+i_0)<2\quad \forall\,\alpha.
$$
Since we also have $i+j-e_\alpha\in\triangle_{d,\Z}$, we get
$$
-2<x_\alpha-2(i_\alpha+j_\alpha+i_0+j_0).
$$
These imply
$$
j_\alpha+j_0<2,\quad \forall\, \alpha.
$$
Summing over $\alpha$ implies $j_0<2$. So $j_0\le 1$. On the other
hand, $d=i_0+j_0\le d-1+j_0$ implies $j_0\ge 1$. We must have
$j_0=1$. But this in turns implies $j_\alpha<2-j_0=1$. So
$j_\alpha\le 0$ and they cannot sum up to get
$\displaystyle\sum_{\beta=1}^{n+1}j_\beta=j_0=1$. This contradiction implies
$i_0=d$ and hence $j_0=0$.

Consider those $\gamma$ such that $i_\gamma\le 1$. (Such $\gamma$
always exists as $i_0=d$ implies $i\neq 0$.) For any one of these
$\gamma$ and any $\alpha\neq\gamma$,
$i+e_\alpha-e_\gamma\in\triangle_{d,\Z}$. So we have
$$
x_\alpha-x_\gamma-2(i_\alpha-i_\gamma) <2.
$$
Using $i+j+e_\gamma-e_\alpha\in\triangle_{d,\Z}$ for all
$\gamma$ and $\alpha$, we have
$$
|x_\gamma-x_\alpha-2(i_\gamma-i_\alpha)-2(j_\gamma-j_\alpha)|<2.
$$
Therefore,
$$
j_\alpha-j_\gamma<2.
$$
Applying this inequality to $\gamma_1$ and $\gamma_2$ with
$i_{\gamma_1}$, $i_{\gamma_2}\ge 1$, we obtain
$$
|j_{\gamma_1}-j_{\gamma_2}|<2.
$$
This implies $j_\gamma\ge 0$ for all $\gamma$ with $i_\gamma\le 1$
or $j_\gamma\le 0$ for all $\gamma$ with $i_\gamma\le 1$.

If $j_\gamma\ge 0$ for all $\gamma$ with $i_\gamma\le 1$, then
$$
0=j_0=\sum_{\{\gamma\,:\,i_\gamma\ge 1\}}j_\gamma
+\sum_{\{\alpha\,:\,i_\alpha=0\}}j_\alpha \ge
\sum_{\{\alpha\,:\,i_\alpha=0\}}1,
$$
as $j_\alpha=i_\alpha+j_\alpha\ge 1$ for $i_\alpha=0$. Hence the set
$\{\alpha\,:\,i_\alpha=0\}$ is empty. So for all $\gamma$,
$i_\gamma\ge 1$ which implies $j_\gamma\ge 0$. Together with
$j_0=0$, we conclude that $j=0$ which is a contradiction.

So we must have $j_\gamma\le 0$ for all $\gamma$ with $i_\gamma\le
1$. Then $j_\alpha-j_\gamma<2$ implies
$$
j_\alpha<2,\quad \forall\, \alpha \mbox{ with }i_\alpha=0.
$$
Therefore
$$
j_\alpha\le 1,\quad \forall\,\alpha \mbox{ with }i_\alpha=0.
$$
On the other hand, for these $\alpha$,
$j_\alpha=i_\alpha+j_\alpha\ge 1$. Hence $j_{\alpha}=1$ for all
these $\alpha$. Putting this back into the inequality, we have
$$
j_\gamma> -1,\quad \forall\,\gamma \mbox{ with }i_\gamma\ge 1.
$$
Hence,
$$
j_\gamma= 0,\quad \forall\,\gamma \mbox{ with }i_\gamma\ge 1.
$$
as $j_\gamma\le 0$ for these $\gamma$. Using $j_0=0$, we conclude
that the set $\{\alpha\,:\,i_\alpha=0\}$ is empty and obtained a
contradiction again. And this completes the proof of Step 4.

\noindent{\bf Step 5:} There exists $\beta$ such that both $i$ and
$i+j$ belong to
$\{y\in\R^{n+1}:y_\beta=0\}\cap\partial\triangle_d\cap\Z^{n+1}$.

\noindent{\em Proof of Step 5:} By Step 4, there exist $\beta$ and
$\alpha$ such that $i_\beta=0$ and $i_\alpha+j_\alpha=0$. If we can
choose $\beta=\alpha$ then we are done. If not, then we have
$$
i_\beta+j_\beta\ge 1\quad\mbox{and}\quad i_\alpha\ge 1.
$$
These imply $i+j-e_\beta$, $i+j+e_\alpha-e_\beta$, $i-e_\alpha$, and
$i+e_\beta-e_\alpha\in\triangle_{d,\Z}$ and hence we have
$$
\begin{cases}
-2<x_\beta-2(i_\beta+j_\beta+i_0+j_0)& \\
x_\alpha-x_\beta-2(i_\alpha-i_\beta)-2(j_\alpha-j_\beta) < 2 & \\
-2 < x_\alpha-2(i_\alpha+i_0) & \\
x_\beta-x_\alpha-2(i_\beta-i_\alpha) < 2
\end{cases}
$$
Therefore
$$
j_\beta-j_\alpha < 2.
$$
Using $i_\beta=0$, $i_\beta+j_\beta\le 1$, and
$i_\alpha+j_\alpha=0$, one get
$$
1+i_\alpha \le j_\beta-j_\alpha < 2,
$$
and arrive at the contradiction that $i_\alpha=0$. This proves the
Step 5.

\noindent{\bf Completion of the proof of the lemma:} By Step 5, if
there exists $x\in\R^{n+1}$ satisfying the condition of the lemma,
then $i$ and $i+j$ belong to $\{y\in\R^{n+1}\,:\,
y_\beta=0\}\cap\triangle_{d,\Z}$ for some $\beta$. This
reduces the argument to one lower dimension. Since the proposition
is clearly true for 1-dimension, induction implies the lemma holds.
\end{proof}

\begin{lem}\label{lem-v2}
For any $i\in\triangle_{d,\Z}$, there exists at most $n+1$
elements $j_{\gamma}\in\{\pm e_\beta, e_\beta-e_\alpha\}_{\beta\neq
\alpha}$ with $j_{\gamma_1}+j_{\gamma_2}\neq 0$ such that there
exists $x\in\R^{n+1}$ satisfying
\begin{eqnarray*}
\langle i, x \rangle -v(i) &=& \langle i+j_{\gamma}, x \rangle
-v(i+j_{\gamma}) \quad \forall,\gamma\\
 & > & \langle r, x \rangle -v(r),\quad r\in \triangle_{d,\Z}\setminus\{i,i+j_{\gamma}\}.
\end{eqnarray*}
\end{lem}
\begin{proof}
We first claim that for any $i\in\triangle_{d,\Z}$ and $j\neq
0 \in\Z^{n+1}$, there exists no $x\in\R^{n+1}$ such that
$$
\langle i, x \rangle -v(i) = \langle i+j, x \rangle -v(i+j) =\langle
i-j, x \rangle -v(i-j).
$$
In fact, if such $x$ exists, then we have the equality
$$
v(i+j)-v(i)=v(i)-v(i-j).
$$
This implies $j=0$ which is a contradiction.

Secondly, we claim that for any $i\in\triangle_{d,\Z}$, there
exists no $x\in\R^{n+1}$ such that
\begin{eqnarray*}
 \langle i, x \rangle -v(i) &=&
\langle i+e_\beta, x \rangle -v(i+e_\beta) \\
&=&\langle i+ e_\alpha, x \rangle -v(i+ e_\alpha)\\
&=&\langle i+e_\beta-e_\alpha, x \rangle -v(i+e_\beta-e_\alpha),
\end{eqnarray*}
and also no $x\in\R^{n+1}$ such that
\begin{eqnarray*}
 \langle i, x \rangle -v(i) &=&
\langle i+e_\beta, x \rangle -v(i+e_\beta) \\
&=&\langle i- e_\alpha, x \rangle -v(i- e_\alpha)\\
&=&\langle i+e_\beta-e_\alpha, x \rangle -v(i+e_\beta-e_\alpha),
\end{eqnarray*}
The first set of equalities implies
$$
\begin{cases}
x_\beta &= 2 i_\beta+2i_0+2 \\
 x_\alpha &= 2 i_\alpha+2i_0 + 2\\
 x_\beta  - x_\alpha   &= 2( i_\beta-i_\alpha )+2 ,
\end{cases}
$$
and the second set implies
$$
\begin{cases}
x_\beta &= 2 i_\beta+2i_0+2 \\
 x_\alpha &= 2 i_\alpha+2i_0 - 2\\
 x_\beta  - x_\alpha   &= 2( i_\beta-i_\alpha )+2 ,
\end{cases}
$$
Both are impossible.

By the two claims, we see that if $\pm e_\beta$ is one of the
$j_{\gamma}$, the $\mp e_\beta$ will not appear in the set
$\{j_{\gamma}\}$; and if $\pm e_\beta$ and $\pm e_\alpha$ belong to
the set $\{j_{\gamma}\}$, then $\pm (e_\beta-e_\alpha)$ will not
appear in the set $\{j_{\gamma}\}$. Therefore, each $\beta=1,\ldots,
n+1$ can appeared once in the set $\{j_{\gamma}\}$ and this
completes the proof of the lemma.

\end{proof}

\noindent{\bf Proof of the key lemma \ref{lem-v}:} It is clear from
the lemmas \ref{lem-v1} and \ref{lem-v2}, the balanced polyhedral
complex $\Pi_v$ corresponding to
$v(j)=\displaystyle\sum_{\beta=0}^{n+1}j_\beta^2$ is a maximal dual complex of
$\triangle_d$ which gives (1) of the lemma. To see (2), we observe
that lemma \ref{lem-v1} implies that any simplex of the subdivision
with an interior point $j$ of the translated simplex
$i-{\i}+\triangle_{n+2}$ as a vertex, then all other vertices belong
to $i-{\i}+\triangle_{n+2}$. Therefore, the subdivision restrict to
a subdivision of $i-{\i}+\triangle_{n+2}$. 

????Finally by lemma \ref{lem-v1}, for each top dimensional face $\frak{F}(i,j)$ given by $i\neq j\in\triangle_{d,\Z}$, we must have $i=j\pm e_\alpha$ or $i=j+e_\alpha -e_\beta$ for some $\alpha$, $\beta=1,\ldots,n+1$. Therefore, if both $i$ or $j\in\partial\triangle_d$, we have $\frak{F}(i,j)$ is unbounded.???? 

This proves the statement (3) and the proof
of the key lemma is completed.


\section{Appendix: Definition of balanced polyhedral complex}
In this appendix, we state the Mikhalkin's definition \cite{Mi03} of
a balanced polyhedral complex for reader's reference.
\begin{defn} A subset $\Pi \in \R^{n+1}$ is called a {\em rational
polyhedral complex} if it can be represented as a finite union of
closed convex polyhedra (possibly semi-infinite) called {\em cells}
in $\R^{n+1}$ satisfying
\begin{enumerate}
\item The slope of the affine span of each cell is rational.
\item If the dimension of the cell is defined to be the dimension
of its affine span and a $k$-dimensional cell is called a {\em
$k$-cell}. Then the boundary of a $k$-cell is a union of
$(k-1)$-cells.
\item Different open cells do not intersect.
\end{enumerate}
\end{defn}

\begin{defn}
\begin{enumerate}
\item The maximum of the dimensions of the cells of a polyhedral
complex $\Pi$ is called the {\em dimension} of $\Pi$. And $\Pi$ is
called a polyhedral $n$-complex if the dimension of $\Pi$ is $n$.
\item A polyhedral $n$-complex is called {\em weighted} if there is a {\em
weight} $w(F)\in \N$ assigned to each of its $n$-cell $F$.
\end{enumerate}
\end{defn}

\begin{defn} For each $n$-cell $F$ of a weighted polyhedral $n$-complex in
$\R^{n+1}$ and an co-orientation on $F$, an integer covector
$$
c_F:\Z^{n+1}\to \Z
$$
is defined by the following conditions
\begin{enumerate}

\item The kernel of $c_F$ is parallel to $F$.

\item The normalized covector $\frac{1}{w(F)}c_F$ is a primitive
integer covector.

\item The covector $c_F$ compatible with the co-orientation of $F$.
\end{enumerate}
\end{defn}

\begin{defn} A weighted polyhedral $n$-complex in $\R^{n+1}$ is called {\em
balanced} if for all $(n-1)$-cell $G\subset \Pi$,
$$
\sum_{s} c_{F_s}=0,
$$
where $F_s$ are the $n$-cells adjacent to $G$ with co-orientation
given by a choice of a rotational direction about $G$.
\end{defn}


\begin{thebibliography}{99}
\bibitem{CL} Chan,; Leung, Conan N. C., {\em Mirror symmetry for toric Fano manifolds via SYZ
transformation}, arXiv: 08012830.

\bibitem{Cohn}  Cohn, P. M., {\em Puiseux's thereom revisited}, J.
Pure Appl. Algebra {\bf 31} (1984), 1-4; correction, {\bf 52}
(1988), 197-198.

\bibitem{EKL} Einsiedler, M.; Kapranov, M.; Lind, D., {\em
Non-Archimedean amoebas and tropical varieties}, J. Reine Angew.
Math. {\bf 601} (2006), 139-157.

\bibitem{Fulton} Fulton, W., {\em Introduction to toric varieties}, Annals of Mathematics Studies, {\bf 131}. The William H. Roever Lectures in Geometry. Princeton University Press, Princeton, NJ, 1993. xii+157pp. ISBN:0-691-00049-2.

\bibitem{Gelfand} I. M. Gelfand, M. M. Kapranov, A. V. Zelevinsky,
{\em Discriminants, resultants, and multidimensional determiniants},
Mathematics: Theory \& Applications. Birkh\"auser Boston, Inc.,
Boston, MA, 1994.

\bibitem{Grif}  Griffiths, P.; Harris, J., {\em Principles of algebraic geometry},Pure and Applied Mathematics. Wiley-Interscience [John Wiley \& Sons], New York, 1978. xii+813pp. ISBN:0-471-32792-1.

\bibitem{Gross} Gross, M., {\em Topological mirror symmetry}, Invent. Math. {\bf 144} (2001), no. 1, 75--137.

\bibitem{HarLaw} Harvey, R.; Lawson, H. B., {\em Calibrated geometries}, Acta Math. {\bf 148} (1982), 47-157.


\bibitem{Mi03} Mikhalkin, G., {\em Decomposition into pairs-of-pants
for complex algebraic hypersurfaces}, Topology {\bf 43} (2004)
1035-1065.

\bibitem{MiNote} Mikhalkin, G.,  {\em Tropical Geometry and Amoebas}, Preprint, 2003.

\bibitem{Ruan} Ruan, W. D., {\em Lagrangian torus fibration of quintic hypersurfaces. I. Fermat quintic case}, Winter School on Mirror Symmetry, Vector Bundles and Lagrangian Submanifolds (Cambridge, MA, 1999), 297--332, AMS/IP Stud. Adv. Math., 23, Amer. Math. Soc., Providence, RI, 2001

\bibitem{SYZ} Strominger, A.; Yau, S-T.; Zaslow, E., {\em Mirror symmetry is $T$-duality}, Winter School on Mirror Symmetry, Vector Bundles and Lagrangian Submanifolds (Cambridge, MA, 1999), 333--347, AMS/IP Stud. Adv. Math., 23, Amer. Math. Soc., Providence, RI, 2001.

\bibitem{Viro}  Viro, O. Ya., {\em Real plane algebraic curves:
constructions with controlled topology}, Leninggrad Math. J. {\bf 1}
(1990), no. 5, 1059-1134.

\bibitem{Yau} Yau, S-T., {\em On the Ricci curvature of a compact K\"ahler manifold and the complex Monge-Amp\`ere equation. I}, Comm. Pure Appl. Math. {\bf 31} (1978), no. 3, 339-411.
\end{thebibliography}
\end{document}